\documentclass[a4paper,11pt,reqno,noindent]{amsart}
\usepackage[centertags]{amsmath}
\usepackage{amsfonts,amssymb,amsthm,dsfont,cases,amscd,esint,enumerate,fancyhdr,newlfont}
\usepackage{comment}
\usepackage[english]{babel}
\usepackage[body={15cm,21.5cm},centering]{geometry} 
\usepackage{tikz}
\usepackage{blindtext}
\usepackage{mdframed}
\usepackage{float}
\usepackage{appendix}
\usepackage{csquotes}
\usepackage{mathtools}
\usetikzlibrary{matrix}
\usepackage{hyperref}
\usetikzlibrary{positioning}
\usepackage{thmtools, thm-restate}
\usetikzlibrary{arrows.meta}

\newtheorem{theor0}{Theorem}[section]

\newenvironment{theor}
{\pushQED{\qed}\begin{theor0}}
	{\popQED\end{theor0}}
\newtheorem{lem0}[theor0]{Lemma}

\newtheorem{prop0}[theor0]{Proposition}
\newenvironment{prop}
{\pushQED{\qed}\begin{prop0}}
	{\popQED\end{prop0}}
\newtheorem{cor0}[theor0]{Corollary}

\newtheorem{propr0}[theor0]{Property}

\newtheorem{hyp0}[theor0]{Hypothesis}

\newtheorem{result0}[theor0]{Result}

\newtheorem{conj0}[theor0]{Conjecture}

\newtheorem{heur0}[theor0]{Heuristics}

\theoremstyle{definition}
\newtheorem{defin0}[theor0]{Definition}
\newenvironment{defin}
{\pushQED{\qed}\begin{defin0}}
	{\popQED\end{defin0}}
\newtheorem{rems0}[theor0]{Remarks}

\newtheorem{ex0}[theor0]{Example}

\newtheorem{exs0}[theor0]{Examples}

\newtheorem{rem0}[theor0]{Remark}
\newenvironment{rem}
{\pushQED{\qed}\begin{rem0}}
	{\popQED\end{rem0}}
\newtheorem{qu0}[theor0]{Question}

\newtheorem{qus0}[theor0]{Questions}

\newtheorem{as0}[theor0]{Assumptions}

\theoremstyle{plain}

\numberwithin{equation}{section}

\newcommand{\N}{\mathbb N}

\newcommand{\Z}{\mathbb Z}
\newcommand{\e}{\varepsilon}
\newcommand{\dist}{\operatorname{dist}}
\newcommand{\R}{\mathbb R}

\newcommand{\calN}{\mathcal N}

\newcommand{\Id}{\operatorname{Id}}

\newcommand{\ul}{{\operatorname{uloc}}}

\newcommand{\Lip}{{\operatorname{Lip}}}
\newcommand{\expec}[1]{\big\langle #1 \big\rangle}

\newcommand{\step}[1]{\noindent \textit{Step} #1.}

\newcommand{\dps}{\displaystyle}

\DeclareMathOperator\supp{supp}
\DeclareMathOperator{\sech}{sech}

\title[Homogenization and shallow water limit of water waves]{Consistency analysis for combined homogenization and shallow water limit of water waves}

\author[A. Gloria]{Antoine Gloria}
\address[Antoine Gloria]{Sorbonne Universit\'e, CNRS, Universit\'e de Paris, Laboratoire Jacques-Louis Lions, 75005~Paris, France \& Universit\'e Libre de Bruxelles, D\'epartement de Math\'ematique, 1050~Brussels, Belgium}
\email{antoine.gloria@sorbonne-universite.fr}

\author[D. Lee]{David Lee}
\address[David Lee]{School of Mathematical and Physical Sciences, University of Technology Sydney (UTS), Sydney, Australia}
\email{davidchanwoo.lee@uts.edu.au}

\begin{document}
	
\begin{abstract}
We consider a shallow water model in a homogenization framework. For periodic topographies, 
Craig, Lannes and Sulem have established a consistency result under some non-resonance conditions. In the present contribution, we significantly relax the periodicity condition and treat general topographies under minimal assumptions.
\end{abstract}
	
	\maketitle

\tableofcontents

\section{Introduction and main result}
	
\subsection{General setting}

We are interested in the motion of an incompressible, irrotational fluid under the influence of gravity, occupying a domain bordered below by a fixed bottom and above by a free surface. These equations are relevant for ocean-coastal dynamics. While the flat-bottom setting is well understood, much less is known when the bathymetry varies. The complicated nature of these equations has motivated the derivation of simpler models which are more tractable and still relevant under specific regimes, and we refer the reader to the extensive monograph \cite{MR3060183} by Lannes on the subject.
	
In this article, we focus on the justification of a shallow water model in a homogenization framework.
Because there are several parameters in the problem (say, to simplify the discussion, the height of the fluid -- the shallow water parameter --, and the period of the oscillating bottom -- the homogenization parameter in the long wave regime -- and its magnitude), there is a whole hierarchy of such models.
In particular, one can first consider the shallow water limit and then the homogenization limit, as done in  \cite{MR693716,MR2430666, MR2121939}.
In \cite{MR2901196}, which is the main inspiration for the present work, Craig, Lannes and Sulem  rather consider a diagonal regime where the shallow water limit and the homogenization limit are taken simultaneously, and interact nontrivially. They perform a consistency analysis (in the wording of \cite{MR3060183}), which amounts to showing that a suitable two-scale ansatz is an approximate solution to the original problem (that is, the ansatz satisfies the equation up to a remainder which is small in the energy norm). This is not yet a complete result since the remainder is not proved to be small in a norm for which the (nonlinear hyperbolic) operator would be continuous, see  \cite[page 30]{MR3060183}. This is a major open problem in the field.
Nevertheless, consistency analysis already reveals some limitations of the derivation when resonance occurs. A resonance is roughly speaking an uncontrolled exchange of energy across scales that potentially destroys the assumption of scale separation made in the very form of the ansatz -- this resonance is interpreted in \cite{MR2901196}  as a nonlinear generalization of the \emph{Bragg resonance}. 
	
The main assumption on the oscillating bottom in \cite{MR2901196} is periodicity. The aim of the present work is twofold. First we revisit the approach of \cite{MR2901196} in order to explain the origin of the ansatz, which indeed combines two arguments: a two-scale expansion (as standard in homogenization) and a small ellipticity contrast limit (the coefficients get closer to the identity as the period vanishes). The first argument allows one to describe oscillations using correctors, whereas the second allows one to obtain a very precise description of them (by linearization). Based on this insight, we turn to the second and main objective of this article: relax the periodicity assumption on the bottom. In a nutshell, we make essentially no assumption on the bottom besides uniform boundedness and an averaging property, and our main result has the flavor of a homogenization result for an algebra with mean value (albeit in the consistency sense, cf.~\cite{MR3060183}), cf.~\cite[Section~7.5]{JKO94}. Due to the small ellipticity contrast limit, the \emph{effective model} is the same whatever the precise description of the bottom. The oscillations do, however, depend on the bottom, and are recovered using (linearized) correctors, which are themselves well-defined in this general setting due to the small ellipticity contrast limit again.
Last, we also extend the crucial non-resonance condition of \cite{MR2901196} to the case of a general bottom. From the analysis point of view, this requires to work with uniformly local Sobolev spaces, rather than on the torus, and raises some subtleties on Fourier multipliers and distributions.

\subsection{Formulation of the water wave problem}

Let $d=1$ (bidimensional problem) or $d=2$ (tridimensional problem). Following \cite[Chapter 1]{MR3060183}, we start with kinematic considerations: the fluid evolves in a domain (an infinite strip) with a free boundary on top and an impermeable wall at the bottom 
\begin{equation}
\Omega(b,\zeta(\cdot,t)):=\{(x,z)\in \mathbb{R}^d\times \mathbb{R}:-1+b(x)<z<\zeta(x,t)\},
\end{equation}
where $\zeta(\cdot,t)$ represents the surface elevation at time $t$ (which is an unknown of the problem) and $b$ represents the variation of the bottom relative to the reference depth $-1$ (which is a datum of the problem). 
We denote by $n$ the unit exterior normal vector on $\partial \Omega(b,\zeta(\cdot,t))$. 
For convenience, we define the fluid height function $h$ via
\begin{equation}\label{eqn:height_function}
h(t,x)=1+\zeta(t,x)-b(x), \quad (t,x)\in [0,T]\times \mathbb{R}^d,
\end{equation}
on some nontrivial time window $[0,T]$.
	
\medskip

We now turn to the physics, and consider the free surface (incompressible and irrotational) Euler equations in $\Omega(b,\zeta(\cdot,t))$ ($\Omega$ in short)
for the fluid velocity $u:\Omega\rightarrow \mathbb{R}^{d+1}$ 
\begin{equation}\label{EQN:EULER}
\begin{cases}
\partial_t u + (u \cdot \nabla)u + \nabla P = -g e_z, \\ \text{div} \, u =0,   \\
\text{curl}\, u =0, 
\end{cases}
\end{equation}
where $P$ is the pressure and $g$ is the gravitational acceleration. 
We complete this system with 
\begin{enumerate}
\item a kinematic condition on the free surface:
\begin{equation*}
\partial_{t}\zeta(t,x)=\sqrt{1+|\nabla \zeta|^2}\, u(t,x,\zeta(t,x)) \cdot n,
\end{equation*}
\item constant atmospheric pressure on the free surface
\begin{equation*}
P(t,x,\zeta(t,x))=P_{const}, 
\end{equation*}
\item impermeable fixed bottom 
\begin{equation} \label{eqn:solid_wall}
u(t,x,-1+b(x))\cdot n=0. 
\end{equation}
\end{enumerate}
We now follow Zakharov/Craig-Sulem to reformulate this system into the water waves problem -- see \cite{zakharov1968stability} and \cite{MR1239970}. 
Since the motion is incompressible and irrotational,  at each time there exists a potential $\Phi:  \Omega\rightarrow \mathbb{R}$ such that $u=\nabla\Phi$. Setting $\psi(t,x):=\Phi(t,x,\zeta(t,x))$ the trace of the velocity potential on the free surface, \eqref{EQN:EULER} and \eqref{eqn:solid_wall} yield the following elliptic boundary value problem at each time $t$
\begin{equation}\label{EQN:ELLIPTIC_PROBLEM}
\left\{
\begin{array}{rcl}
\triangle_{d+1} \Phi(t,\cdot)&=&0,\quad \text{in }\Omega(b,\zeta(\cdot,t)),\\
\Phi(t,x,\zeta(t,x))&=&\psi(t,x),
\\
 \partial_{n}\Phi(t,x,-1+b(x))&=&0,
\end{array}
\right.
\end{equation}
where $\triangle_{d+1}$ is the Laplacian in $\R^{d+1}$.
This elliptic problem is then used to define the Dirichlet-to-Neumann operator $G[\zeta,b]$ acting on functions $\psi$  (at each time $t$) through $\Phi$
via 
\begin{equation*}
G[\zeta,b]\psi (t,x):=\sqrt{1+|\nabla \zeta|^2}\, \partial_n\Phi(t,x,\zeta(t,x)),
\end{equation*}
where $\nabla$ denotes the gradient in $\R^d$.
By the chain rule, the above system and its boundary conditions reduce to the water waves problem:
\begin{equation}\label{EQN:water_waves}
\begin{cases}
&\partial_{t}\zeta-G[\zeta,b]\psi=0,\\
&\partial_t\psi+g\zeta+\frac{1}{2}|\nabla \psi|^2-\frac{(G[\zeta,b]\psi+\nabla\zeta \cdot \nabla \psi)^2}{2(1+|\nabla \zeta|^2)}=0,
\end{cases}
\end{equation}
where $g$ is the gravitational constant.
System \eqref{EQN:water_waves} models water waves above an oscillating bottom, and we still need to choose the regime we shall consider. There are essentially four parameters in the problem:
\begin{enumerate}
\item the typical amplitude $A$ of surface waves,
\item the typical wave length $\lambda$ of surface waves,
\item the typical amplitude $B$ of the variations of the bottom $b$,
\item the typical wavelength  $\ell$ of the bottom. 
\end{enumerate}
Using the dimensionless variables  
\begin{equation*}
x=\lambda x', \quad z= z', \quad t=\frac{\lambda}{\sqrt{g}}t',\quad \zeta=A\zeta',\quad \Phi= A\lambda \sqrt{g}\Phi', \quad b= {B}b'(x/\ell),
\end{equation*}
we obtain the four independent dimensionless parameters  
\begin{equation}
\begin{split}
&\mu=\lambda^{-2}\quad \text{(shallowness parameter)}, \quad \varepsilon=A \quad\text{(nonlinearity parameter)},\\
&\beta=B \quad \text{(topography parameter)},\quad\, \gamma=\frac{\ell}{\lambda} \quad\text{(transversality parameter)}. 
\end{split}
\end{equation}
In dimensionless form and with $g=1$, the system for \eqref{EQN:water_waves} becomes 	
\begin{equation}\label{EQN:water_waves_dimensionless}
\begin{cases}
&\partial_{t}\zeta-\frac{1}{\mu}G_{\mu}[\zeta,\beta b(\frac\cdot\gamma)]\psi=0\\
&\partial_t\psi+\zeta+\frac{1}{2}|\nabla \psi|^2-
\mu\frac{(\frac{1}{\mu}G_{\mu}[\zeta,\beta b(\frac\cdot\gamma)]\psi+\nabla\zeta \cdot \nabla \psi)^2}{2(1+\mu|\nabla \zeta|^2)}=0,
\end{cases}
\end{equation}
where
\begin{equation*}
G_{\mu}[\zeta,\beta b(\tfrac\cdot\gamma)]\psi=\sqrt{1+|\nabla \zeta|^2}  \partial_n^\mu \Phi(\cdot,\zeta),
\end{equation*}
where $\Phi$ now solves the following elliptic problem
\begin{equation}\label{EQN:ELLIPTIC_PROBLEM_dimensionless}
\left\{
\begin{array}{rcl}
(\mu \triangle+ \partial_z^2) \Phi &=&0,\quad \text{on }\Omega,\\
\Phi(\cdot,\zeta)&=&\psi, \\
\partial_n^\mu \Phi(\cdot,-1+\beta b(\tfrac\cdot\gamma))&=&0
\end{array}
\right.
\end{equation}
in the fluid domain $\Omega=\Omega(b,\zeta):=\{(x,z)\in \R^{d+1}\,|\,-1+\beta b(\frac x\gamma)<z<\zeta(x)\}$,
where $\triangle$ denotes the Laplacian in $\R^d$, $n$ denotes the outward normal, $\partial_n^\mu:=n \cdot \nabla^\mu$, and  
\[
\nabla^\mu := \begin{pmatrix}
\mu \nabla \\
\partial_z
\end{pmatrix}.
\]
It is with this formulation \eqref{EQN:water_waves_dimensionless} \&~\eqref{EQN:ELLIPTIC_PROBLEM_dimensionless} of the water waves problem that we will work in this article.

\subsection{Main result}

Following \cite{MR2901196}, we investigate the regime $\beta=\sqrt{\mu}=\gamma\ll1$ with no smallness of the amplitude of surface waves, that is, $\e=1$. The fact that $\beta=\gamma$ corresponds to small bathymetry slope, while the roughness strength is
$\rho:=\frac{\sqrt\mu}\gamma=1$. The water waves system thus simplifies to finding $(\zeta,\psi):[0,T]\times \R^d \to \R\times \R$ that satisfy
\begin{equation}\label{e4} 
\begin{cases}
&\partial_{t}\zeta-\frac{1}{\mu}G_{\mu}[\zeta,\sqrt \mu b(\frac \cdot{\sqrt \mu})]\psi=0\\
&\partial_t\psi+\zeta+\frac{1}{2}|\nabla \psi|^2-
\mu\frac{(\frac{1}{\mu}G_{\mu}[\zeta,\sqrt \mu b(\frac \cdot{\sqrt \mu})]\psi+\nabla\zeta \cdot \nabla \psi)^2}{2(1+\mu|\nabla \zeta|^2)}=0,
\end{cases}
\end{equation}
where
\begin{equation}\label{e5} 
G_{\mu}[\zeta,\sqrt \mu b(\tfrac\cdot{\sqrt \mu})]\psi=\sqrt{1+|\nabla \zeta|^2}   \partial_n^\mu \Phi(\cdot,\zeta),
\end{equation}
and $\Phi$ solves the  elliptic problem
\begin{equation}\label{e6}  
\left\{
\begin{array}{rcl}
(\mu \triangle+ \partial_z^2)\Phi&=&0,\quad \text{on }\Omega,\\
\Phi(\cdot,\zeta)&=&\psi, \\
 \partial_n^\mu \Phi(\cdot,-1+\sqrt \mu b(\tfrac\cdot{\sqrt \mu}))&=&0
\end{array}
\right.
\end{equation}
in the fluid domain $\Omega=\Omega(b,\zeta):=\{(x,z)\in \R^{d+1}\,|\,-1+\sqrt \mu b(\tfrac x{\sqrt \mu})<z<\zeta(x)\}$.

\medskip

Although it is quite clear, let us emphasize that the above system is a coupled system of two nonlinear nonlocal equations.
Without the shallow water limit, such a system would be extremely difficult (if not impossible) to homogenize. In particular, one would not expect a separation of scales.
What makes the above system easier to handle is the combined shallow water and homogenization limits. Indeed, due to the shallow water regime, the large-scale behavior of the solution becomes independent of the oscillations of the bottom. This forces separation of scales, and in turn allows us to characterize correctors and reconstruct the oscillations of the solution. Of course, things are made complicated by the nonlinearity of the system, which imposes nontrivial conditions on the large-scale behavior of the solution to prevent constructive (and uncontrolled) interferences with the oscillating bottom (which is called nonlinear Bragg resonance in \cite{MR2901196}). Last, the result is not a full homogenization result because of the poor existence and regularity theory available for the original system.
 
\medskip

Let us start with a two-scale expansion as used in homogenization theory,
and denote by $(t,x) \mapsto (\zeta_0,\psi_0)(t,x)$ the ``expected'' large-scale behavior of $(\zeta,\psi)$, and set $V_0:=\nabla \psi_0$ and $h_0:=1+\zeta_0$. 
If there is separation of scales, one  expects $(\zeta,\psi)$ to be close to two-scale expansions of the form 
\begin{equation}\label{e3}
\begin{cases}
&\zeta^{2s}_{\mu}(t,x)=\zeta_0(t,x)+\sqrt{\mu}\zeta_c(\frac{x}{\sqrt{\mu}},h_0(t,x),V_0(t,x)),\\ 
&\psi^{2s}_{\mu}(t,x)=\psi_0(t,x)+\mu \psi_c(\frac{x}{\sqrt{\mu}},h_0(t,x),V_0(t,x)),
\end{cases}
\end{equation}
where the correctors $\zeta_c$ and $\psi_c$ are  there to reconstruct oscillations due to the bottom $b$. In \eqref{e3}, we made
two assumptions on these correctors: they add only \textit{spatial} oscillations at scale $\sqrt{\mu}$ (that is, there are no time oscillations in this ansatz -- which is a  limitation), and the large-scale solution $(\zeta_0,V_0)$ at point $x$ and time $t$ enters as a parameter. In fact, for all values  $(\bar h,\bar V)$ (that is, an element of $\R\times \R^d$) of the parameters, $\zeta_c(\cdot,\bar h,\bar V)$ and $\psi_c(\cdot,\bar h,\bar V)$ should solve a system of equations depending on $b$ and on these parameters.

\medskip

Since $\lim_{\mu\downarrow 0} \sqrt{\mu}b(\frac \cdot {\sqrt \mu}) =0$, it is natural to expect that $(\zeta_0,\psi_0)$ should satisfy the shallow water wave system with flat bottom, which -- see~Section~\ref{sec:motiv} -- takes the form
\begin{equation}\label{e8}
\begin{cases}
\hspace{1.2cm}\partial_{t}\zeta_0+\nabla\cdot (V_0 h_0)&=0,
\\
\partial_{t}V_0+\nabla\zeta_0+(V_0\cdot \nabla)V_0 &=0,
\end{cases}
\end{equation}
where $h_0:=1+\zeta_0$ is the effective height function and $V_0=\nabla \psi_0$, with suitable initial conditions.
Assume that the solution exists (and is continuous) on some time interval $[0,T]$.

\medskip

An informal analysis -- see~Section~\ref{sec:motiv} -- will then lead us to an explicit system for the correctors  $\zeta_c$ and $\psi_c$, which is given by 
\begin{equation}\label{eqn:corrector_eqns}
\begin{pmatrix}
\bar V\cdot\nabla_y & -|\nabla_y|\tanh{(\bar h|\nabla_y|)}\\
I & \bar V\cdot\nabla_y
\end{pmatrix}
\begin{pmatrix}
\zeta_c\\
\psi_c
\end{pmatrix}	=\begin{pmatrix}
(\bar V\cdot\nabla_y) \text{sech}(\bar h|\nabla_y|)b\\
0
\end{pmatrix}, 
\end{equation}
for $(\bar h,\bar V)\in X_{T}{\supset}\{(h_{0}(t,\cdot),V_0(t,\cdot)), t \in [0,T]\},$  which we can assume to be an open set (since $V_0$ and $h_0$ will be continuous).
Inserting the second line into the first line of \eqref{eqn:corrector_eqns} yields
\begin{equation}\label{eqn:psi_corrector}
-\left ((\bar V\cdot \nabla_y)^2+|\nabla_{y}|\tanh{(\bar h|\nabla_y|)}\right)\psi_c=(\bar V\cdot \nabla_y) \sech(\bar h|\nabla_y|)b. 
\end{equation}
Once $\psi_c$ is constructed through \eqref{eqn:psi_corrector}, $\zeta_c$ is immediately recovered from $\zeta_c=-(V_0\cdot \nabla_y)\psi_c$. One of our main contributions concerns the solvability of the corrector equations for general oscillating bottoms $b$, for which we now introduce some assumptions.

By a formal application of the Fourier transform, it will become apparent that one needs to impose assumptions on $(b,X_{T})$ such that
$\xi \mapsto \frac{\sech(\bar h|\xi|)}{(\bar V\cdot \xi)^2-|\xi|\tanh(\bar h|\xi|)}(i\xi)$ makes sense on a neighborhood of the support of $\hat b$ for all $(\bar V,\bar h)\in X_{T}$.

\begin{defin}\label{def:non-resonance}
Let $\hat{b}$ denote the Fourier transform of $b$ in the distributional sense, and decompose $\supp \hat b:=\Xi_1\cup \Xi_2$ for some discrete set $\Xi_1$.
We say that $(b, X_{T})$ satisfies the non-resonance condition of order $(m,n) \in \N^2_0$ if 
there exist a continuous function $\hat{K}_{0}:\R^d\times X_{T}\to \R^d$ and an open set $\Xi_3$ containing $\Xi_2$ such that 
for all $\xi \in \Xi_1 \cup \Xi_3$ and all $(\bar V,\bar h)\in X_T$ we have
\[
\hat{K}_{0}(\xi,\bar V,\bar h)\,=\,\frac{\sech(\bar h|\xi|)}{(\bar V\cdot \xi)^2-|\xi|\tanh(\bar h|\xi|)}(i\xi),
\]
and that satisfies the following bounds, together with its derivatives.
For  $\alpha \in \N^d_0$, set\footnote{with the understanding that
$(i\xi)^\alpha=\prod_{j=1}^d (i\xi_j)^{\alpha_j}$.}
 $\hat{K}_{\alpha}(\xi,\bar V,\bar h):=(i\xi)^\alpha\hat{K}_{0}(\xi,\bar V,\bar h)$. There exists an $\e>0$ such that for all $(\xi,\bar V,\bar h)\in \R^d\times X_{T}$
\begin{equation}\label{eqn:prop_8_resonance}
|\partial_{\xi}^\tau\partial_{\bar{V},\bar{h}}^\beta\hat{K}_{\alpha}(\xi,\bar V,\bar h)|\leq 
\frac{C_{\alpha,\beta, \tau}}{(1+|\xi|)^{d+\e}},
\end{equation}
where $|\tau|\leq d+1$, $|\alpha|\leq m,|\beta|\leq n$ and  where $C_{\alpha,\beta, \tau}$ is a finite constant independent of $\xi\in \R^d$ and $(\bar V,\bar h)\in X_{T}$.
We take the convention of denoting 
\begin{equation}
\partial_{\xi}^\tau\partial_{\bar{V},\bar h}^\beta\hat{K}_{\alpha}:=
\partial_{\xi}^\tau
\partial_{\bar V}^{\beta_{1}}
\partial_{\bar h}^{\beta_2}\hat{K}_{\alpha}, \quad \text{for $\tau \in \mathbb{N}_0^d, \beta=(\beta_1,\beta_2)\in \mathbb{N}_0^{d+1}$.}
\end{equation}
\end{defin}
In essence, the non-resonance condition is a requirement that enables us to identify $\hat{K}_{\alpha}$ as a relatively \enquote{nice} pseudo-differential symbol, like the  H\"ormander class of symbols~\cite{MR2304165}, cf.~\cite[Section 7]{MR3465379} for a setting involving uniformly local Sobolev spaces. This motivates the introduction of the sets $\Xi_1$ and $\Xi_2$. For the discrete set $\Xi_1$, the structure theorem for zero-th order distributions imposes little restrictions. For the non-discrete set $\Xi_2$ however, it does not suffice to know the multiplier on $\Xi_2$ to characterize the distribution, whence the condition on the larger set $\Xi_3$.

\begin{restatable}[]{as0}{nonresonance}\label{assumption:main_thm}
We assume that 
\begin{itemize}
\item the oscillating bottom $b:\R^d \to \R$ satisfies $b \in W^{1,\infty}(\R^d)$ and $0\notin \supp \hat b$,
\item $(V_0,\zeta_0)\in L^\infty([0,T];H^s(\mathbb{R}^d)^{d+1})\cap {\Lip}([0,T];H^{s-1}(\mathbb{R}^d)^{d+1})$ for some $s>\frac{d}{2}+2$.
\item $(b,X_{T})$ satisfies the non-resonance condition of order $(2,2)$, see Definition \ref{def:non-resonance},
\item For some $\alpha_0>0$ we have that $h_0=1+\zeta_0\geq \alpha_0$ in $[0,T]$. 
\end{itemize}
\end{restatable}
Under Assumptions \ref{assumption:main_thm}, we are able to construct correctors $(\zeta_c,\psi_c)$, see Theorems~\ref{thm:corrector_existence} and~\ref{thm:corrector_sobolev}, which leads us to the following consistency result.
\begin{theor}\label{thm:main}
Let $(\zeta_0,\psi_0):[0,T] \times \R^d\to \R\times \R$ be a smooth solution of \eqref{e8} decaying at infinity on 
some time interval $[0,T]$. Under Assumptions~\ref{assumption:main_thm},  the two-scale ansatz~\eqref{e3} is well-defined as a function of class $H^2(\R^d)$ and it is an approximate solution of~\eqref{e4}--\eqref{e6} in the sense that 
\begin{equation*} 
\begin{cases}
&\partial_{t}\zeta_{\mu}^{2s}-\frac{1}{\mu}G_{\mu}[\zeta_{\mu}^{2s},\sqrt \mu b(\frac \cdot{\sqrt \mu})]\psi_{\mu}^{2s}=E^1_\mu, \\
&\partial_t\psi_{\mu}^{2s}+\zeta_{\mu}^{2s}+\frac{1}{2}|\nabla \psi_{\mu}^{2s}|^2-
\mu\frac{(\frac{1}{\mu}G_{\mu}[\zeta_{\mu}^{2s},\sqrt \mu b(\frac \cdot{\sqrt \mu})]\psi_{\mu}^{2s}+\nabla\zeta_{\mu}^{2s} \cdot \nabla \psi_{\mu}^{2s})^2}{2(1+\mu|\nabla \zeta_{\mu}^{2s}|^2)}=E^2_\mu,
\end{cases}
\end{equation*}
where
\begin{equation*}
\sup_{t\in [0,T]} \|E^{1}_{\mu} \|_{L^2(\R^d)}\lesssim \mu^\frac{3}{8}, \quad \sup_{t\in [0,T]} \|E^{2}_{\mu}\|_{H^\frac{1}{2}(\R^d)}\lesssim  \mu^\frac{3}{4}.
\end{equation*}
\end{theor}
In a nutshell, provided the large-scale solution $(\zeta_0,\psi_0)$ does not resonate with the oscillating bottom $b$, one can construct a two-scale ansatz that is an approximate solution of the original Zakharov/Craig-Sulem system. Recent Bloch--Floquet work (for periodic $b$) also
confirms the dynamical relevance on band gaps of the non-resonance condition for the linearized operator with periodic bottom oscillations \cite{LacaveMenardSulem2025}. 
To complete this result, let us first recall existence and regularity results for~\eqref{e8}.
\begin{theor}\cite[Proposition 6.1]{MR3060183} 
Let $r_0>\frac{d}{2}$, $s\geq r_0+1$ and $(\zeta^\circ,V^\circ)\in H^{s}(\mathbb{R}^d)\times H^{s}(\mathbb{R}^d)^{d}$. Assume that $1+\zeta^\circ\geq \alpha_0$ for some $\alpha_0>0$. Then there exists a time $T>0$ such that~\eqref{e8} has a unique solution $(V_0,\zeta_0)$ with initial condition $(\zeta^\circ,V^\circ)$. Moreover, 
we have that $(V_0,\zeta_0)\in L^\infty([0,T];H^s(\mathbb{R}^d)^{d+1})\cap {\Lip}([0,T];H^{s-1}(\mathbb{R}^d)^{d+1})$.  
\end{theor}
Second, the error term in Theorem~\ref{thm:main} is controlled in the energy norm of the system: Theorem~\ref{thm:main}  is a consistency result, not a convergence theorem.  To turn the estimates
\[
 \|E^1_\mu\|_{L^\infty(0,T;L^2)}\lesssim \mu^{3/8},
 \qquad
 \|E^2_\mu\|_{L^\infty(0,T;H^{1/2})}\lesssim \mu^{3/4}
\]
into a smallness result for the exact solution with zero initial data, one would need a forced stability estimate for the full water-wave system, uniform in the oscillating bottom
$B_\mu(x):=\sqrt\mu\,b(x/\sqrt\mu)$.

For fixed sufficiently smooth bottoms, the local Cauchy theory and the uniform shallow-water/asymptotic theory yield this type of nonlinear stability statement in high Sobolev energies; see in particular~\cite{Lannes2005,MR3060183,MR2372806}.  These results are uniform with respect to the usual dimensionless parameters only for families of bottoms bounded in the required smooth norms.  They do not apply uniformly here: even if \(b\) is smooth,
\[
 \|\partial^\alpha B_\mu\|_{L^\infty}
 \sim \mu^{(1-|\alpha|)/2}\|\partial^\alpha b\|_{L^\infty},
\]
so every estimate depending on \(W^{k,\infty}\) or high Sobolev norms of the bottom loses uniformity as soon as \(k\geq2\).  Moreover, the remainders in Theorem~\ref{thm:main} are small only in the low natural norms above; differentiating the two-scale ansatz costs powers of \(\mu^{-1/2}\), and the high Sobolev residuals required by standard quasilinear estimates are not small.

The general-bottom Cauchy theories of Alazard--Burq--Zuily and the Hadamard continuous dependence result of Nguyen are robust with respect to rough or non-localized domains \cite{AlazardBurqZuily2014,MR3465379,Nguyen2016}, but they are not formulated as forced shallow-water stability estimates uniform for the above oscillating family of bottoms. 
 
The small-data two-dimensional series of works of Hunter--Ifrim--Tataru, Ifrim--Tataru, and Ai--Ifrim--Tataru gives sharp modified-energy, normal-form, and global small-data results for flat infinite-depth gravity waves~\cite{HunterIfrimTataru2016,IfrimTataru2016,AiIfrimTataru2026,AiIfrimTataru2022}. The finite-depth flat-bottom result of Harrop-Griffiths--Ifrim--Tataru proves local well-posedness and cubic lifespan bounds for small data, uniformly in the infinite-depth limit \cite{HarropGriffithsIfrimTataru2017}; the Morawetz estimate of Alazard--Ifrim--Tataru is likewise uniform in the infinite-depth limit, under the assumption of scale-invariant bounds~\cite{AlazardIfrimTataru2022}.  These works seem to support strong small-data stability mechanisms, but they do not treat non-flat oscillatory bathymetry, and their uniformity is not with respect to \(B_\mu\), as would be needed here.

\subsection{Discussion of the non-resonance condition and examples}

Theorem~\ref{thm:main} is only valid under the non-resonance condition of Assumptions~\ref{assumption:main_thm}, which is quite implicit. In this paragraph, we present several sufficient conditions for its validity, on several concrete examples of oscillating topographies.
As a first sanity check, we quickly show that the non-resonance condition of Craig, Lannes and Sulem in \cite{MR2901196} in the periodic setting indeed implies the validity of Assumptions~\ref{assumption:main_thm}. Then, we investigate which topographies can be considered in practice beyond the periodic setting of \cite{MR2901196}. Can we for instance consider standard models of random topography (an open problem mentioned in \cite[Section 1.7.2]{MR3060183})? For both questions, we take inspiration from the discussion of \cite[Section~4.3]{MR2901196}.

\medskip

In what follows, we denote by $b$ the topography, $\hat b$ its Fourier transform (as a distribution), and 
set $\Xi:=\supp{\hat b}$. We also recall that we assume the lower
bound $1+\zeta_0\geq \alpha_0$ for some $\alpha_0>0$.

\medskip

We start with the comparison to the  $(0,1)^d$-periodic setting considered in \cite{MR2901196}, 
for which $\Xi \subset 2\pi \Z^d \setminus \{0\}$. The non-resonance condition of \cite{MR2901196} takes the following guise: There exists $0<\alpha_*< \alpha_0$ such that we have for all $\xi \in \Xi$ and $(\bar{h},\bar{V})\in X_{T}$,
\begin{equation}\label{e.res-cond-Lannes}
	|(\bar V\cdot \xi)^2-|\xi|\tanh(\bar h|\xi|)|>e^{-\alpha_{*}|\xi|}.
\end{equation}
Let us check that condition~\eqref{e.res-cond-Lannes} implies the validity of Assumptions~\ref{assumption:main_thm}. 
If $b$ is $(0,1)^d$-periodic, then for all $\xi \in \Xi$, $\dist(\xi,\Xi \setminus\{\xi\}) \ge 2\pi$.
We shall consider the more general case when $\Xi=\Xi_1$ is discrete (that is, $\Xi_2=\varnothing$), bounded away from $0$, and satisfies 
\begin{equation}\label{e.dioph}
\dist(\xi,\Xi \setminus\{\xi\}) \ge c e^{-\frac{\alpha_0-\alpha_{*}}{d+2}|\xi|},
\end{equation}
for some $c >0$. This covers for instance all quasi-periodic topographies, and even some almost-periodic topographies, two cases which are not addressed in \cite{MR2901196}.
It remains to construct $\hat K_0$ and prove \eqref{eqn:prop_8_resonance}:
\begin{itemize}
\item Let $\rho_*:\R^d \to [0,1]$ be a smooth bump function with support in $B_1=B(0,1)$ and such that $\rho_*(0)= 1$;
\item For all $\xi \in \Xi$, define $r_*(\xi)=1 \wedge \dist(\xi,\Xi \setminus\{\xi\})$, and set $\rho_\xi:=\rho_*(\frac2{r_*(\xi)} (\cdot-\xi))$;
\item Define
\[
\hat{K}_{0}(\xi,\bar V,\bar h)\,:=\,\sum_{\xi'\in \Xi} \frac{\sech(\bar h|\xi'|)}{(\bar V\cdot \xi')^2-|\xi'|\tanh(\bar h|\xi'|)}(-i \xi)\rho_{\xi'}(\xi).
\]
\end{itemize}
By definition, for all $\xi\in \Xi$, 
\[
\hat{K}_{0}(\xi,\bar V,\bar h)\,=\frac{\sech(\bar h|\xi|)}{(\bar V\cdot \xi)^2-|\xi|\tanh(\bar h|\xi|)}(-i \xi).
\]
In addition, for all $\tau\le d+1$ and $\beta_1,\beta_2\le 2$,
\[
\partial_{\xi}^\tau\partial_{\bar{V},\bar{h}}^\beta\hat{K}_{\alpha}(\xi,\bar V,\bar h)
\,=\, \sum_{\xi'\in \Xi} \partial_{\bar{V},\bar{h}}^\beta\Big(\frac{\sech(\bar h|\xi'|)}{(\bar V\cdot \xi')^2-|\xi'|\tanh(\bar h|\xi'|)}\Big)
\partial_{\xi}^\tau \Big(-(i \xi)^{\alpha+1}\rho_{\xi'}(\xi)\Big),
\]
and a direct calculation yields
\[
|\partial_{\xi}^\tau\partial_{\bar{V},\bar{h}}^\beta\hat{K}_{\alpha}(\xi,\bar V,\bar h)|
\,\le \, C e^{-\frac{\alpha_0-\alpha_{*}}{d+3}|\xi|},
\]
where $C$ depends on $d$ and on $\sup_{X_T}( |\bar V|+|\bar h|)$, which implies \eqref{eqn:prop_8_resonance}, and thus the validity of Assumptions~\ref{assumption:main_thm}.

\medskip

We then turn to more general topographies.
Let $b_*$ be a non-negative bump function in the Schwartz class, which we assume to be compactly supported away from zero in Fourier space. Set $\Xi_*:=\supp{b_*}$, $r_*:=\inf_{\xi \in \Xi_*} |\xi|$ and $R_*:=\sup_{\xi \in \Xi_*} |\xi|$,  $V_*:=\sup |\bar V|$, and recall that $\bar h\ge h_0>0$. 
Assume that for some $0<\kappa<r_*$ we have
\begin{equation}\label{froude}
\inf_{\rho \in [r_*-\kappa,R_*+\kappa]} \frac{\tanh(h_0 \rho)}{V_*^2 \rho}>1.
\end{equation}
If \eqref{froude} holds, then $(b,X_T)$ satisfies Assumptions~\ref{assumption:main_thm},
for topographies $b$ of the form
\begin{equation*}
b(x)=\sum_{z\in \mathbb{Z}^d}\gamma_{z}b_{*}(x-x_z),
\end{equation*}
for any $\gamma_j \in [0,1]$ and any point set $\{x_z\}_{z\in \Z^d}$ such that $\sup_{z\in \Z^d} \sharp\{j:x_j \in Q(z)\}<\infty$, $Q(z)=z+[-1/2,1/2)^d$.
This covers for instance random amorphous models ($\gamma_j=1$ and $\{x_z\}_{z\in \Z^d}$ a stationary ergodic hardcore point set), the random displacement model ($\gamma_j=1$ and $x_z=z+\tau_z$ where $\tau_z$ are iid translations in $[0,1)^d$), and the random ``Anderson'' model ($\gamma_j$ iid, and $x_z=z$).
In this case, one can choose $\Xi_1=\varnothing$ and $\Xi_2=\Xi \subset \Xi_*$. 
Set $\Xi_3:=\{\xi \in \R^d \,|\, \dist(\xi,\Xi_*)<\kappa\}$ for some $0<\kappa<r_*$ such that \eqref{froude} holds.
We can define
\[
\hat K_0(\xi,\bar V,\bar h)\,:=\, \frac{\sech(\bar h|\xi|)}{(\bar V\cdot \xi)^2-|\xi|\tanh(\bar h|\xi|)}(-i \xi)\rho(\xi),
\]
where $\rho:\R^d\to [0,1]$ is a smooth cut-off function such that $\rho|_{\bar \Xi_*}\equiv 1$, $\rho|_{\R^d \setminus  \Xi_3}\equiv 0$,
and $|\nabla^k \rho|\lesssim \kappa^{-k}$ for all $k\le d+1$. 
Because of \eqref{froude},
\[
\inf_{\xi \in \supp{\rho}, (\bar V,\bar h) \in X_T} |\xi|\tanh(\bar h |\xi|)-(\bar V \cdot \xi)^2 >0.
\]
Combined with the fact that $\supp{\rho}$ is bounded, this implies \eqref{eqn:prop_8_resonance}, and thus 
Assumptions~\ref{assumption:main_thm}.
Of course, \eqref{froude} is a strong assumption that amounts to small data (for the limiting equation), and indeed prevents any Bragg resonance. Condition~\eqref{froude}  is somehow related to a lower bound on the Froude number $Fr:= \frac{\bar V^2}{\bar h}$.

\subsection{Outline of article}

The rest of this contribution is organized as follows.
In Section~\ref{sec:motiv} we motivate the specific form \eqref{e3} of the ansatz, of the equations for $(\zeta_0,\psi_0)$, and of the correctors.
In Section~\ref{sec:corrector} we address the well-posedness of the corrector equations, and finally 
prove Theorem~\ref{thm:main} in Section~\ref{sec:proof}.

\subsection*{Notation}
\begin{itemize}
\item We ley $x\in \mathbb{R}^d$, for $d=1,2$, denote the horizontal variables and $z$ denote the vertical variable.
\item For a given multi-index $\alpha\in \mathbb{N}^{d}_0$ we will denote $|\alpha|=\sum_{k=1}^d\alpha_k$ and $\partial_{x}^\alpha=\prod_{k=1}^{d}\partial_{x_k}^{\alpha_{k}}$.  We also set $\langle x\rangle:=(1+|x|^2)^\frac{1}{2}$ for $x\in \mathbb{R}^d$. 
\item We take the following convention for the Fourier transform and its inverse:
\begin{equation*}
\mathcal{F}(\varphi)(\xi):=\hat{\varphi}(\xi):=\int_{\mathbb{R}^d}e^{-ix\,\cdot \xi}\varphi(x)\,dx,\quad \text{for $\varphi \in \mathcal{S}(\mathbb{R}^d)$ and $\xi\in \mathbb{R}^d$}, 
\end{equation*}
where $\mathcal{S}(\mathbb{R}^d)$ is the set of \emph{Schwartz functions}, and
\begin{equation*}
\varphi(x)=\mathcal{F}^{-1}(\hat{\varphi})(x):=(2\pi)^{-d}\int_{\mathbb{R}^d}\hat{\varphi}(\xi)e^{ix\cdot \xi}\,d\xi,\quad \text{for $\varphi \in \mathcal{S}(\mathbb{R}^d)$ and $x\in \mathbb{R}^d$}. 
\end{equation*}
\item A similar convention will be used on $\mathcal{S}'(\mathbb{R}^d)$ which is the space of tempered distributions. Specifically, for $u\in \mathcal{S}'(\mathbb{R}^d)$, we define the Fourier transform of $u$ as 
\begin{equation*}
\langle \hat{u},\varphi\rangle_{\mathcal{S}',\mathcal{S}}:=\langle u,\hat{\varphi}\rangle_{\mathcal{S}',\mathcal{S}}, \quad \text{for $\varphi\in \mathcal{S}(\mathbb{R}^d)$}. 
\end{equation*}
\item We denote $H^{s}(\mathbb{R}^d):=\{u\in \mathcal{S}'(\mathbb{R}^d):\|u\|_{H^s}<\infty\},$ where  
\begin{equation*}
\|u\|_{H^s}:=\|(1-\Delta)^{s/2}u\|_{L^2},
\end{equation*}
and $\| .\|_{L^p}$ is the usual $L^p(\mathbb{R}^d)$ norm for $p\in [1,\infty]$. 
\item We use the standard notation $f(\nabla)$ for the Fourier multiplier having the symbol $f(-i\xi)$.  
\item Often we will consider \emph{multi-scale functions}. Specifically, they are of the form $f(x,y)$, for which the realization (for a given $\gamma \in (0,1)$) is given as $f(x,x/\gamma)$, or $f(x,x/\gamma,z)$. The variable $y\in \mathbb{R}^d$ is known as the \emph{fast variable}. We add the subscript $x$ or $y$ to be specific as to what we are differentiating with respect to. For instance, 
	\begin{equation*}
	(\nabla f) (x,x/\gamma)= (\nabla_{x}f)(x,x/\gamma)+\tfrac{1}{\gamma}(\nabla_{y}f)(x,x/\gamma), \quad \text{for $x\in \mathbb{R}^d$.}
	\end{equation*}
\item We will also use the short-hand notation 
	\begin{equation}\label{e.not-nab-mu}
		\nabla^\mu_{x,z}=(\sqrt{\mu}\nabla_{x},\partial_{z})^{T},\quad \nabla^\mu_{y,z}=(\sqrt{\mu}\nabla_{y},\partial_{z})^{T}, \quad \nabla^\mu_{x,0}=(\sqrt{\mu}\nabla_{x},0)^{T}.
	\end{equation}
When $\mu=1$ we simply set $\nabla^\mu_{x,z}=\nabla_{x,z}$, $\nabla^\mu_{y,z}=\nabla_{y,z}$ and $\nabla^\mu_{x,0}=\nabla_{x,0}$. 
\end{itemize}

\section{Formal asymptotic analysis}\label{sec:motiv}
	
\subsection{General strategy: two-scale expansion of the system}

We first assume the height function $h$ and the fluid velocity $\nabla \psi$ are given by a leading-order term plus some corrections of order $\gamma$ and frequency $\gamma^{-1}$, which leads to the ansatz
\begin{equation}\label{eqn:multi-scale-ansatz_2}
\begin{cases}
&\zeta^{2s}_{\gamma}(t,x)=\zeta_0(t,x)+\gamma \zeta_1(t,x,\tfrac x\gamma) ,\\ 
&\psi^{2s}_{\gamma}(t,x)=\psi_0(t,x)+\gamma^2 \psi_1(t,x,\tfrac x\gamma),     
\end{cases}
\end{equation}
where $\zeta_1$ and $\psi_1$ are multiscale functions that have the joint averaging property with $b$ (that is, any local function of $\psi_1$, $\zeta_1$, and $b$ still has the averaging property\footnote{In particular, if $b$ is periodic (resp. random stationary ergodic), then $\psi_1$ and $\zeta_1$ are periodic (resp. random stationary ergodic) in the fast variable -- however we shall see that such strong structural assumptions are not needed.}), and  we say that a multiscale function $h$ has ensemble average $\expec{h}$ if the following limit exists for all $t,x$:
\[
\expec{h(t,x,\cdot)}:=\lim_{R \to \infty} \fint_{[-R,R)^d} h(t,x,y)dy
\]
(here, $\expec{\cdot}$ denotes an ensemble average, not the Japanese bracket).
On the one hand, this implies that $x\mapsto \zeta_1(t,x,\tfrac x \gamma)$ and $x\mapsto \psi_1(t,x,\tfrac x\gamma)$ converge weakly to some limits as $\gamma \downarrow 0$. 
We shall assume that, as corrections, $\expec{\psi_1}\equiv 0$ and $\expec{\zeta_1}\equiv 0$.
On the other hand, this also implies that $\expec{\nabla_y h} \equiv 0$ (provided that $\nabla_y$ and $\expec{\cdot}$ commute).
We first insert this ansatz into the water wave system, and sort the terms of the equations in function of their powers of the small parameter $\gamma$. 
We first simplify the system using that $\frac{1}{\mu}G_{\mu}[\zeta^{2s}_\gamma,\sqrt \mu b(\frac \cdot{\sqrt \mu})]\psi^{2s}_\gamma$ is of order 1, and obtain with the choice $\sqrt\mu=\gamma \ll 1$ for the first equation\footnote{The loose notation $o$ (resp. $O$) means that identities hold true up to that order $o$ (resp. $O$) in all the norms we shall need (in particular, we shall use that this holds at the level of gradients $\nabla^\mu$).}
\begin{equation}\label{g1}
\partial_{t}\zeta_0+\tfrac{1}{\mu}G_{\mu}[\zeta^{2s}_\gamma,\sqrt \mu b(\tfrac \cdot{\sqrt \mu})]\psi^{2s}_\gamma=O(\sqrt \mu).
\end{equation}
The asymptotic analysis of  $\frac{1}{\mu}G_{\mu}[\zeta^{2s}_\gamma,\sqrt \mu b(\frac \cdot{\sqrt \mu})]\psi^{2s}_\gamma$ will result in an expression involving $\zeta_1$ and $\psi_1$. In order to close the system, we need the second equation to involve $\zeta_1$ and $\psi_1$ as well, which requires us to go one order further in the accuracy, in the form of 
\begin{equation}\label{g2}
\partial_t\psi_0+(\zeta_0+\gamma \zeta_1)+\frac{1}{2}|\nabla \psi_0|^2+\gamma \nabla \psi_0 \cdot \nabla_y \psi_1 =o(\gamma).
\end{equation}
Our main task is now to identify the leading order term of the Dirichlet-to-Neumann map 
$\frac{1}{\mu}G_{\mu}[\zeta^{2s}_\gamma,\sqrt \mu b(\frac \cdot{\sqrt \mu})]\psi^{2s}_\gamma$,
and hopefully close the system.

\subsection{Asymptotic expansion of the Dirichlet-to-Neumann operator}

For the analysis, it is convenient to reformulate the Dirichlet-to-Neumann operator on a fixed domain by a change of variables, and then proceed with the asymptotic expansion.

\subsubsection{Reformulation on a fixed domain}

Set $\Omega_0:=\{(x,z)\in \R^{d+1}\,|\, x \in \R^s, z\in (-1,0)\}$ and 
$\Omega_\gamma(t):=\{(x,z)\in \R^{d+1}\,|\, x\in \R^d, \ -1+\gamma b(\frac x {\gamma} )<z<\zeta^{2s}_\gamma(t,x)\}$.
In this paragraph, since the time $t$ is a parameter we shall drop it from our notation.
Following \cite{MR3060183}, provided $\zeta^{2s}_{\gamma},b\in W^{1,\infty}(\mathbb{R}^d)$ we can transform the elliptic problem \eqref{EQN:ELLIPTIC_PROBLEM_dimensionless} on $\Omega_\gamma$ into an elliptic problem on the flat strip $\Omega_0$ by using the diffeomorphism
\begin{equation*}
\begin{array}{rcrcl}
S_\gamma^{2s}&:&\Omega_0&\rightarrow& \Omega_\gamma\\
&&(x,z)&\mapsto &(x,z+\sigma^{2s}_\gamma(x,z)),
\end{array}
\end{equation*}
where 
\begin{equation}\label{eqn:sigma}
\sigma^{2s}_\gamma(x,z):=\underbrace{(1+z)\zeta_0(x)}_{\dps =:\sigma_0(x,z)}+ \gamma \underbrace{\Big( (1+z)\zeta_1(x,\tfrac x\gamma)-   zb(\tfrac x \gamma)\Big)}_{\dps =:\sigma_1(x,z,y)|_{y=\frac x \gamma}}.
\end{equation}
In what follows, all multiscale functions are evaluated at $(x,z,\frac{x}{\gamma})$.
The function $\phi:=\Phi\circ S_\gamma^{2s}$ now solves 
\begin{equation}\label{eqn:elliptic_strip}
\left\{
\begin{array}{rcl}    
\nabla^\mu \cdot P[\sigma_\gamma^{2s}]\nabla^\mu \phi&=&0, \quad \text{on $\Omega_0$},\\
\phi|_{z=0}&=&\psi^{2s}_\gamma,\\
 e_z\cdot P[\sigma^{2s}_{\gamma}]\nabla^\mu\phi|_{z=-1}&=&0,
\end{array}
\right.
\end{equation}
where $P[\sigma^{2s}_{\gamma}]$ is the nonlinear function of $\sigma^{2s}_\gamma$
\begin{equation*}
P[\sigma^{2s}_{\gamma}]=\left(\begin{array}{cc}
h^{2s}_{\gamma} \Id_d & -\sqrt{\mu} \nabla \sigma_{\gamma}^{2s}\\
-\sqrt{\mu} (\nabla \sigma_{\gamma}^{2s})^T&
\frac{1 +\mu\vert \nabla \sigma_{\gamma}^{2s}\vert^2}{h^{2s}_{\gamma}}
\end{array}\right),
\end{equation*}
and  the height function $h^{2s}_{\gamma}$ is given by
\begin{equation*}
h^{2s}_{\gamma}(x):=1+\partial_{z}\sigma_{\gamma}^{2s}(x)=\underbrace{1+\zeta_0(x)}_{\dps =: h_0(x)}
+ \gamma \underbrace{\Big( \zeta_1(x,\tfrac x\gamma)-   b(\tfrac x \gamma)\Big)}_{\dps =:h_1(x,y)|_{y=\frac x \gamma}}.
\end{equation*}
(Recall that $\nabla$ denotes the gradient in $\R^d$, whereas $\nabla^\mu$ is the gradient in $\R^{d+1}$ scaled by $\sqrt{\mu}$ in the $x$-coordinate.)
All in all, the Dirichlet-to-Neumann map now takes the form 	
\begin{equation}\label{DTN:transformed}
G[\zeta_{\gamma}^{2s},\gamma b(\tfrac \cdot \gamma)]\psi^{2s}_\gamma=e_{z}\cdot P[\sigma_{\gamma}^{2s}]\nabla^\mu\phi|_{z=0}.
\end{equation}

As explained above, we shall decompose $\phi$ into a part due to the shallow-water expansion and a part due to the oscillating bottom. To this aim we decompose  $P[\sigma^{2s}_{\gamma}]$ as
\begin{eqnarray}
P[\sigma^{2s}_{\gamma}]&=&P_0+\sqrt{\mu}P_1 \label{a.P}
\\
&=&P_{0,0}+\sqrt{\mu} (P_{0,1}+P_{1,0}) + \mu P_{1,1},\label{a.Pb}
\end{eqnarray}
where $P_0=P_{0,0}+\sqrt{\mu} P_{0,1}$, $P_1= P_{1,0}+\sqrt \mu P_{1,1}$, and 
\[
P_{0,0}:=\begin{pmatrix}
h_0\Id_d & 0\\
0 & \frac{1}{h_0}
\end{pmatrix},
\quad P_{0,1} :=\begin{pmatrix}
0\Id_d & -(1+z)\nabla\zeta_0\\
-(1+z)(\nabla\zeta_0)^T & \sqrt{\mu}\frac{|\nabla \zeta_0|^2}{h_0}
\end{pmatrix},
\]
and
\[
P_{1,0}:=
\begin{pmatrix}
h_1 \Id_d  & -\nabla_{y} \sigma_1\\
- \nabla_{y} \sigma_1^T &
\frac{b-\zeta_1}{h_0^2}
\end{pmatrix},
\quad P_{1,1}:=\begin{pmatrix}
0 \Id_d  & -\nabla_{x} \sigma_1\\
- \nabla_{x} \sigma_1^T &
p_{1122} 
\end{pmatrix},
\]
with $p_{1122}:=\sqrt{\mu}^{-1} \Big(\sqrt{\mu}^{-1}(\frac{1+\mu |\nabla \sigma^{2s}_\gamma|^2}{h_\gamma^{2s}}-\frac{1+\mu |\nabla \sigma_0|^2}{h_0})-\frac{b-\zeta_1}{h_0^2}\Big)$.
In the analysis $\mu P_{1,1}$ is higher-order and will be considered as an error term.

\medskip

Using \eqref{a.P} and \eqref{eqn:multi-scale-ansatz_2} (and the relation $\gamma^2=\mu$), we decompose $\phi=\phi_0+ \mu \phi_1$, where $\phi_0$ and $\phi_1$ solve
\begin{equation}\label{e:shallow}
\left\{
\begin{array}{rcl}    
\nabla^\mu \cdot P_0 \nabla^\mu \phi_0&=&0 , \quad \text{on $\Omega_0$},\\
\phi_0|_{z=0}&=&\psi_0,\\
 e_z\cdot P_0\nabla^\mu\phi_0|_{z=-1}&=&0
\end{array}
\right.
\end{equation}
and 
\begin{equation}\label{e:hom}
\left\{
\begin{array}{rcl}    
\nabla^\mu \cdot P[\sigma^{2s}_{\gamma}] \nabla^\mu \phi_1&=&-\tfrac1{\sqrt \mu}\nabla^\mu \cdot P_1\nabla^\mu \phi_0 , \quad \text{on $\Omega_0$},\\
\phi_1|_{z=0}&=&\psi_1,\\
 e_z\cdot P[\sigma^{2s}_{\gamma}] \nabla^\mu\phi_1|_{z=-1}&=&-\tfrac1{\sqrt \mu}e_z\cdot P_1\nabla^\mu\phi_0|_{z=-1}.
\end{array}
\right.
\end{equation}
The asymptotic analysis of \eqref{e:shallow} when $\mu \downarrow 0$ is the standard shallow water limit, whereas the asymptotic analysis of \eqref{e:hom} is of the homogenization type in the small ellipticity contrast regime.

\subsubsection{Shallow water limit: asymptotic analysis of~\eqref{e:shallow}}\label{sec:swl}

In this subsection, we display the classical justification (see \cite[Section 3.2]{MR2901196}  and  \cite[Section 3.6]{MR3060183}) that the solution of \eqref{e:shallow} is well-approximated by 
\begin{equation}\label{g.6}
\phi_0\,=\,\psi_0+\mu \phi_{0,1}+O(\mu^2), \quad \phi_{0,1}(x,z):=-(\tfrac{z^2}{2}+z) h_0^2(x) \triangle \psi_0(x).
\end{equation} 
Direct calculations yield
\begin{eqnarray*}
\nabla^\mu \cdot P_0 \nabla^\mu \psi_0&=& \mu h_0 \triangle \psi_0,
\\
\nabla^\mu \cdot P_0\nabla^\mu \phi_{0,1}&=&-h_0 \triangle \psi_0 - \mu \Psi(h_0,\zeta_0,\psi_0),
\end{eqnarray*}
where 
\begin{multline*}
\Psi(h_0,\zeta_0,\psi_0)\,=\,(\tfrac{z^2}2+z) \nabla \cdot (2h_0^2\triangle \psi_0 \nabla \zeta_0 +h_0^3 \nabla \triangle \psi_0) -(z+1)^2 \nabla \cdot (h_0^2 \triangle \psi_0 \nabla \zeta_0) 
\\
+ |\nabla \zeta_0|^2h_0 \triangle \psi_0
- (\tfrac32 z^2+3z+1) (2h_0 |\nabla \zeta_0|^2 \triangle \psi_0+h_0^2 \nabla \zeta_0 \cdot \nabla \triangle \psi_0),
\end{multline*}
and we also have
\[
e_z \cdot P_0 \nabla^\mu \psi_0|_{z=-1}=0, \quad e_z \cdot P_0 \nabla^\mu \phi_{0,1}|_{z=-1}=0.
\]
All in all, $\phi_0-(\psi_0+\mu \phi_{0,1})$ thus satisfies 
\begin{equation*}
\left\{
\begin{array}{rcl}    
\nabla^\mu \cdot P_0 \nabla^\mu (\phi_0-\psi_0-\mu \phi_{0,1})&=&\mu^2\Psi(h_0,\zeta_0,\psi_0), \quad \text{on $\Omega_0$},\\
 (\phi_0-\psi_0-\mu \phi_{0,1})|_{z=0}&=&0,\\
 e_z\cdot P_0\nabla^\mu  (\phi_0-\psi_0-\mu \phi_{0,1})|_{z=-1}&=&0,
\end{array}
\right.
\end{equation*}
from which \eqref{g.6} follows (provided $\Psi(h_0,\zeta_0,\psi_0)$ is smooth enough).

\subsubsection{Homogenization and asymptotic analysis of~\eqref{e:hom}}\label{sec:hom}

We argue, following \cite{MR2901196}, that 
\begin{equation}\label{g.7}
\phi_1=\phi_{1,0}+ O(\sqrt \mu),
\quad
\left\{
\begin{array}{rcl}    
(h_0^2\Delta_y+\partial_{z}^2) \phi_{1,0}&=&0 , \quad \text{on $\Omega_0$},\\
\phi_{1,0}|_{z=0}&=&\psi_1,\\
\frac{1}{h_0}\partial_{z}\phi_{1,0}|_{z=-1} &=&\nabla_{y} b  \cdot \nabla \psi_0.
\end{array}
\right.
\end{equation}
To this aim, we rewrite \eqref{e:hom} as
\begin{equation}\label{a1} 
\left\{
\begin{array}{rcl}    
\nabla^\mu \cdot P_{0,0} \nabla^\mu \phi_1&=&R_1, \quad \text{on $\Omega_0$},\\
\phi_1|_{z=0}&=&\psi_1,\\
 e_z\cdot P_{0,0} \nabla^\mu\phi_1|_{z=-1}&=&- e_z\cdot P_{1,0}\nabla^\mu\psi_0|_{z=-1}+R_2,
\end{array}
\right.
\end{equation}
where
\begin{eqnarray*}
R_{1}&:=&- \nabla^\mu \cdot P_{1,1}\nabla^\mu \psi_0 - \tfrac1{\sqrt \mu}\nabla^\mu \cdot P_{1,0}\nabla^\mu \psi_0-\tfrac1{\sqrt \mu}\nabla^\mu \cdot P_1\nabla^\mu (\phi_0-\psi_0)
\\
&&- \sqrt{\mu} \nabla^\mu\cdot (P_{0,1}+P_1)\nabla^\mu \phi_1,
\\
R_2&:=&- e_z \cdot P_{1,1}\nabla^\mu \psi_0|_{z=-1}
-\tfrac1{\sqrt \mu} e_z\cdot P_1 \nabla^\mu(\phi_0-\psi_0)|_{z=-1}
\\
&&- \sqrt \mu e_z\cdot (P_{0,1}+P_1) \nabla^\mu\phi_1|_{z=-1}.
\end{eqnarray*}
Let us use the following notation: $\nabla_{y,z}:=\begin{pmatrix} \nabla_y \\ \partial_z \end{pmatrix}$, that is the gradient in the fast horizontal variables and the usual gradient in the vertical variable, and $\nabla_{x,0}:= \begin{pmatrix} \nabla_x \\ 0 \end{pmatrix}$, that is the gradient in the slow horizontal variables and zero in the vertical variable. This notation allows us to decompose the scaled gradient $\nabla^\mu$ of a multiscale function $(x,z,y) \mapsto \eta_1(x,z,y)$ evaluated at $(x,z,\frac x{\sqrt \mu})$ as
\begin{equation}\label{aa0}
\nabla^\mu \eta_1 \,=\, \nabla_{y,z} \eta_1+\sqrt{\mu} \nabla_{x,0} \eta_1,
\end{equation}
whereas for functions $\eta_0$ that depend on the slow variable $x$ but not on the vertical variable $z$ nor on the fast variables $y$, one has
\begin{equation}\label{aa1}
\nabla^\mu \eta_0 \,=\,  \sqrt{\mu} \nabla_{x,0} \eta_0.
\end{equation}
With these explicit scalings in $\mu$ and \eqref{g.6}, it is easy to see that the remainders $R_1$ and $R_2$ are expected to be small and can thus be neglected (to leading order) in \eqref{a1}, and that one can replace $\nabla^\mu$ by $\nabla_{y,z}$ (to leading order) in the left-hand side of  \eqref{a1} -- all this with accuracy $\sqrt{\mu}$.
Noting that by  \eqref{g13} (with $z=-1$ in which case $- \nabla_y  \sigma_1=\nabla_y b$ by \eqref{eqn:sigma}), one may reformulate the Neumann boundary condition in \eqref{a1} as
\[
-e_z\cdot P_{1,0}\nabla^\mu\psi_0|_{z=-1}\,=\, \nabla_y b\cdot \nabla \psi_0,
\]
we have derived~\eqref{g.7}.

\subsubsection{Two-scale expansion of the Dirichlet-to-Neumann map}\label{sec:D2N}
In this paragraph, we argue that 
\begin{equation}\label{eqn:DTN_ref}
\tfrac{1}{\mu}G_{\mu}[\zeta_{\mu}^{2s},\sqrt{\mu} b(\tfrac \cdot {\sqrt{\mu}})]\psi_{\mu}^{2s}=G_{\text{eff}}+\sqrt{\mu}\,G_{\text{res}},
\end{equation}
where $G_{\text{res}}=O(1)$ and 
\begin{equation}\label{eqn:DTN:eff0}
G_{\text{eff}}:=-\nabla\cdot(h_0V_0)-(V_0\cdot \nabla_{y})\zeta_1+\frac{1}{h_0}\partial_{z}\phi_{1,0}|_{z=0}.
\end{equation}
The starting point is the combination of~\eqref{g.6} and~\eqref{g.7}, which yields the approximation 
\begin{equation}\label{g.8}
\phi(x,z)=\psi_0(x)+\mu(\phi_{0,1}(x,z)+\phi_{1,0}(x,z,\tfrac x\gamma))+O(\mu^{3/2}).
\end{equation}
Taking the scaled gradient of \eqref{g.8}, that is, $\nabla^\mu \phi=\nabla^\mu \psi_0+\mu(\nabla^\mu \phi_{0,1}+\nabla^\mu \phi_{1,0})+O(\mu^{3/2})$, we obtain using \eqref{a.Pb}, 
\begin{eqnarray*}
\lefteqn{\tfrac{1}{\mu} G[\zeta_{\mu}^{2s},\sqrt{\mu}b(\tfrac \cdot {\sqrt{\mu}})]\psi_{\mu}^{2s}}
\\
&=&
\tfrac{1}{\mu}e_{z}\cdot P[\sigma_{\mu}^{2s}]\nabla^\mu\phi|_{z=0}
\\
&=&\tfrac{1}{\mu}e_{z}\cdot P[\sigma_{\mu}^{2s}]\nabla^\mu \Big(\psi_0+\mu(\phi_{0,1}+\phi_{1,0})\Big)|_{z=0}+O(\sqrt \mu)
\\
&=&\tfrac{1}{\mu}e_{z}\cdot P[\sigma_{\mu}^{2s}]\nabla^\mu \psi_0|_{z=0}
+e_{z}\cdot P_{0,0} \nabla^\mu \phi_{0,1}|_{z=0}
+e_{z}\cdot P_{0,0} \nabla^\mu\phi_{1,0}|_{z=0}+O(\sqrt \mu).
\end{eqnarray*}
The first two right-hand side terms will combine into the standard Dirichlet-to-Neumann map in the shallow water limit, whereas the third term yields the contribution of the oscillating bottom.
Let us start with the latter. Since 
$P_{0,0} =\begin{pmatrix}
h_0\Id_d & 0\\
0 & \frac{1}{h_0}
\end{pmatrix}$,
we obtain the last right-hand side term of \eqref{eqn:DTN:eff0}
\begin{equation}\label{g10}
e_{z}\cdot P_{0,0} \nabla^\mu\phi_{1,0}|_{z=0}=\frac{1}{h_0}\partial_{z}\phi_{1,0}|_{z=0}.
\end{equation}
Using \eqref{g.6}, the second right-hand side term yields
\begin{equation}\label{g.11}
e_{z}\cdot P_{0,0} \nabla^\mu\phi_{0,1}|_{z=0}=\frac{1}{h_0}\partial_{z}\phi_{0,1}|_{z=0}=-h_0 \triangle \psi_0.
\end{equation}
It remains to analyze the contribution of the first term. Again we use the decomposition~\eqref{a.Pb} of $P[\sigma_\gamma^{2s}]$, and claim that only $\sqrt\mu (P_{0,1}+P_{1,0})$ contributes to leading order.
Indeed, since $\psi_0$ does not depend on $z$ and $P_{0,0}$ is diagonal, $e_z \cdot P_{0,0} \nabla^\mu \psi_0 \equiv 0$, whereas $e_z \cdot P_{1,1} \nabla^\mu \psi_0 = - \sqrt \mu \nabla_x \sigma_1 \cdot \nabla \psi_0$ so that 
$\frac1\mu e_z \cdot \mu P_{1,1} \nabla^\mu \psi_0 = O(\sqrt \mu)$ is higher-order. We then turn to the leading-order terms. First,
\begin{equation} \label{g12}
\tfrac1\mu e_{z}\cdot \sqrt \mu P_{0,1} \nabla^\mu\psi_0 |_{z=0}\,=\, -\nabla \zeta_0 \cdot \nabla \psi_0=-\nabla h_0 \cdot \nabla \psi_0.
\end{equation}
Second, 
\begin{equation} \label{g13}
\tfrac1\mu e_{z}\cdot \sqrt \mu P_{1,0} \nabla^\mu\psi_0 |_{z=0}\,=\, - \nabla_y  \sigma_1 \cdot  \nabla \psi_0
\stackrel{\eqref{eqn:sigma}}=-\nabla_y \zeta_1(x,\tfrac x \gamma) \cdot \nabla \psi_0(x).
\end{equation}
The combination of \eqref{g10}--\eqref{g13} then yields the claimed formula~\eqref{eqn:DTN_ref}.

\subsection{The homogenized system}

We come back to the general strategy, and recall \eqref{g1} and \eqref{g2}
\begin{eqnarray*}
\partial_{t}\zeta_0+\tfrac{1}{\mu}G_{\mu}[\zeta^{2s}_\gamma,\sqrt \mu b(\tfrac \cdot{\sqrt \mu})]\psi^{2s}_\gamma&=&O(\sqrt \mu),
\\
\partial_t\psi_0+(\zeta_0+\gamma \zeta_1)+\frac{1}{2}|\nabla \psi_0|^2+\gamma \nabla \psi_0 \cdot \nabla_y \psi_1 &=&o(\gamma).
\end{eqnarray*}
To conclude the asymptotic analysis, we have to be more precise for the Dirichlet-to-Neumann map \eqref{eqn:DTN_ref} \&~\eqref{eqn:DTN:eff0}.
Since \eqref{g.7} has constant coefficients, it can be solved explicitly in Fourier space (with respect to the fast variable only), and we have the pseudo-differential formula (provided it makes sense)
\begin{equation}\label{g.9}
\phi_{1,0}(x,z,y)=\frac{\cosh(h_0(x)(z+1)|\nabla_y|)}{\cosh(h_0(x)|\nabla_y|)}
\psi_1(x,y) + \frac{\sinh(h_0(x) z|\nabla_y|)}{\cosh(h_0(x)|\nabla_y|)} \frac{V_0(x)\cdot \nabla_{y}}{|\nabla_y|}b(y).
\end{equation}
With \eqref{g.9} we may reformulate $G_{\text{eff}}$ as
\begin{equation}\label{eqn:DTN:eff}
G_{\text{eff}}=-\nabla\cdot(h_0V_0)-(V_0\cdot \nabla_{y})\zeta_1+|\nabla_y| \tanh(h_0\,|\nabla_y|) \psi_1+ (V_0\, \cdot \nabla_y) \sech(h_0\,|\nabla_y|) b.
\end{equation}
All in all, this entails 
\begin{eqnarray}
{O(\sqrt{\mu})}
&=&\partial_{t}\zeta_0-\nabla\cdot(h_0V_0) \nonumber \\
&&-(V_0\cdot \nabla_{y})\zeta_1+|\nabla_y| \tanh(h_0\,|\nabla_y|) \psi_1+ (V_0\, \cdot \nabla_y) \sech(h_0\,|\nabla_y|) b,\label{e.ageff1}
\\
o(\gamma)&=&\partial_t\psi_0+(\zeta_0+\gamma \zeta_1)+\frac{1}{2}|\nabla \psi_0|^2+\gamma \nabla \psi_0 \cdot \nabla_y \psi_1 .\label{e.ageff2}
\end{eqnarray}
We are now in the position to derive an effective system. 
For the leading-order term, we first take the ensemble average $\langle \cdot \rangle$ of \eqref{e.ageff1} \&~\eqref{e.ageff2}, use that $\expec{\zeta_1}=\expec{\psi_1}=\expec{b}=0$, and take the limit $\gamma=\sqrt \mu \downarrow 0$ to the effect that 
\begin{eqnarray*}
\partial_t \zeta_0+\nabla\cdot(h_0V_0)&=&0,\\
\partial_t\psi_0+\zeta_0 +\frac{1}{2}|\nabla \psi_0|^2 &=&0.
\end{eqnarray*}
This entails the effective system \eqref{e8} (upon differentiating the second equation), and it remains to characterize $\zeta_1$ and $\psi_1$.
Using the effective system \eqref{e8}, the equations~\eqref{e.ageff1} \&~\eqref{e.ageff2} take the simpler form, after taking the limit $\gamma=\sqrt \mu \downarrow 0$,
\begin{equation}\label{g4}
-(V_0\cdot \nabla_{y})\zeta_1+|\nabla_y|\tanh(h_0|\nabla_y|)\psi_1+(V_0\cdot \nabla_y) \sech(h_{0}|\nabla_y|)b=0
\end{equation}
and
\begin{equation}\label{g5}
 \zeta_1+ \nabla \psi_0 \cdot \nabla_y \psi_1 =0.
\end{equation}
In these equations, $t,x$ are parameters (through $V_0(t,x)$ and $h_0(t,x)$) and the differential operators are with respect to the fast variable $y$. In particular, 
\begin{equation}\label{g00}
\zeta_1(t,x,y)=\zeta_c(y,h_0(t,x), V_0(t,x)), \quad \psi_1(t,x,y)=\psi_c(y,h_0(t,x), V_0(t,x)),
\end{equation}
where $\zeta_c$ and $\psi_c$ are the correctors given by the system \eqref{eqn:corrector_eqns}.
This motivates the form of Theorem~\ref{thm:main}.

\medskip

A few comments are in order.
First, as one can see on \eqref{g4} and \eqref{g5}, the correctors are nonlinear functions of $h_0$ and $V_0=\nabla \psi_0$ -- this is not surprising in view of the nonlinearity of the system. Second, the correctors solve a constant-coefficient linear system with forcing term depending linearly on $b$. In turn, this implies that we indeed have $\expec{\zeta_1}=\expec{\psi_1}=0$ as we assumed from the very beginning.

\subsection{Relation to classical homogenization}

There are at least three peculiar features about Theorem~\ref{thm:main} and the above informal derivation of the two-scale expansion.

\medskip

First, we aim to replace the hyperbolic system \eqref{e4}--\eqref{e6} by the simpler effective system \eqref{e8}, together with the oscillatory correction \eqref{e3}. In particular, this implies that there is a separation of scales, and that the hyperbolic system homogenizes to a hyperbolic system. In view of the literature (e.g.~\cite{E-92,Jabin-Tzavaras-09,Dalibard-09}) on homogenization of scalar conservation laws, this separation of scales could be surprising. 
 
\medskip

Second, the homogenized system \eqref{e8} is the standard shallow water system on a flat bottom, as if the oscillating bottom had no effect on the limiting equation. This is also rather surprising since in the case of the homogenization of the 2D Euler system \cite{DG-24}, for which in good settings there is a separation of scales, the limiting system is not an Euler system (the vorticity formulation homogenizes but is not associated with an Euler equation any longer). The only exception to this is in the perturbative regime (say, small density of obstacles), in which case, to leading-order, the homogenized system is indeed the Euler system.

\medskip

Third, the result is quantitative: the remainders $E^1_\mu,E^2_\mu$ in the two-scale expansion error are estimated in terms of $\mu$ in Theorem~\ref{thm:main}. This is again surprising since the result is stated without quantitative assumptions on the oscillating bottom $b$ -- only with the qualitative averaging property of Assumptions~\ref{assumption:main_thm} (and a non-resonance condition on $V_0$, $h_0$ and $b$). In general, the quantification of the homogenization error comes from the growth of correctors. In the elliptic setting for instance, correctors are sublinear solutions to equations posed on the whole space, and their actual growth strongly 
depends on the structural assumptions made on the coefficients. To be precise, if the coefficients are periodic, then one can look for periodic correctors (which thus do not grow at infinity), that do exist by Poincar\'e's inequality on the torus. In case of random coefficients with long-range correlations, correctors do indeed grow (sublinearly) at infinity \cite{GNO-reg,GNO-quant}.
In the present setting, although the corrector problem is again posed on an unbounded domain, that domain is a strip and the equation is completed by a Dirichlet boundary condition on one boundary. Hence, we have a Poincar\'e inequality at our disposal regardless of the oscillations of the boundary, which indeed yields bounded correctors, and explains why we have quantitative results without structural assumptions on the oscillating bottom $b$.

\medskip

The reason why the three slightly surprising features above can hold is the perturbative regime due to the choice of parameters in the system \eqref{e4}--\eqref{e6}, and in particular the small amplitude of the bathymetry of the bottom combined with the smallness of the shallowness parameter\footnote{The non-perturbative regime would be an oscillating bottom of the form $b(\frac \cdot \gamma)$, whereas we have here $\gamma b(\frac \cdot \gamma)$.}.
\begin{itemize} 
\item  The corrector equations~\eqref{g4} and~\eqref{g5} are linear and depend linearly on $b$, and they have constant coefficients. The linearity of the equation comes from the shallow-water limit (which linearizes the Dirichlet-to-Neumann operator with respect to correctors). The linearity of the equation with respect to $b$ comes from the small amplitude of the bathymetry -- which amounts to a problem with small ellipticity contrast. Yet, the corrector equation still depends nonlinearly on the value $(\bar V,\bar h)$ of the large-scale solution and it is not solvable for all values of $(\bar V,\bar h)$ -- we call these bad values resonances. Fortunately, since the equation is linear with respect to $b$ and has constant coefficients (because the reference bottom is flat), it can be explicitly solved in Fourier space, and one therefore has a way to characterize resonances.
\item The homogenized operator is in fact independent of $b$ as a consequence of the shallow-water limit.
Hence, it is well-defined for all $(\bar V,\bar h)$. In particular, we can define a good candidate $(V_0,h_0)$ for the large-scale behavior of the solution by solving the effective equation, and check a posteriori whether the latter does  have resonances or not. If not, one can consider the associated correctors, and in turn get the consistency result. 
\end{itemize}
The key to prove the consistency for combined homogenization and shallow water limit of water waves is the 
explicit formulas for correctors in function of $b$ in Fourier space. 
Although this is completely unrelated to our results, note that explicit formulas for correctors in Fourier space are also crucial in \cite{CMOW-25} to study critical random drift-diffusion in $\R^d$. The small parameter is played there by the P\'eclet number.

%%%%%%%%

\section{Construction of the two-scale expansion}\label{sec:corrector}

The goal of this section is the solvability of the corrector problem~\eqref{eqn:corrector_eqns}, which is at the core of our analysis, and constitutes the main originality with respect to~\cite{MR2901196}.
We start by introducing uniformly-local Sobolev spaces, and then turn to the solvability of~\eqref{eqn:corrector_eqns} in these spaces.

\subsection{Uniformly local Sobolev spaces}

Assumptions~\ref{assumption:main_thm} entail $b\in W^{1,\infty}(\R^d)$. Yet, to solve the corrector problem and characterize the Neumann map, we rather consider $b$ as an element of the uniformly local Sobolev space  $H^1_{\ul}(\mathbb{R}^d)$ that we now define. 

\begin{defin}\label{definition:uniformly_local}
Let $\omega\in \mathcal{D}(\mathbb{R}^d)$ where $\text{supp}(\omega)\subset [-1,1]^d$, $\omega=1$  on $\{x\in \R^d:|x|\leq \tfrac{1}{2}\}$ such that $\sum_{q\in \mathbb{Z}^d}\omega_q(x)=1$, for $x\in \mathbb{R}^d$ where $\omega_{q}:=\omega(\cdot -q)$ for $q\in \mathbb{Z}^d$. 
For $s\in \mathbb{R}$, we say that $u\in H^{s}_{\ul}(\mathbb{R}^d)$ if $\omega_{q}u\in H^s(\mathbb{R}^d)$ for all $q\in \mathbb{Z}^d$, and 
\begin{equation*}
\|u\|_{H^s_{\ul}(\mathbb{R}^d)}:=\sup_{q\in \mathbb{Z}^d}\|\omega_qu\|_{H^s(\mathbb{R}^d)}<\infty. 
\end{equation*}
We call $H^{s}_{\ul}(\mathbb{R}^d)$ a uniformly local Sobolev space. For $s=0$ we denote $H^{0}_{\ul}(\mathbb{R}^d)=L^2_{\ul}(\mathbb{R}^d)$ with corresponding norm given as $\|.\|_{L^2_{\ul}(\mathbb{R}^d)}$.
Similarly, we also say that\footnote{Note that the weight $\omega$ only depends on the horizontal variable $x\in \R^d$, not on the vertical variable $z$ which lives in the compact set $[-1,0]$.} $u\in H^1_{\ul}(\Omega_0)$ if 
\begin{equation*}
\|u\|_{H^1_{\ul}(\Omega_0)}:=\sup_{q\in \mathbb{Z}^d}\|\omega_qu\|_{H^1(\Omega_0)}<\infty.
\end{equation*}
We also say that $u\in H_{0,\ul}^{1}(\Omega_0)$ if $u\in H^1_{\ul}(\Omega_0)$ and $\omega_{q}u \in H^1_0(\Omega_0)$ (that is, $\omega_{q}u|_{z=0}\equiv 0$ in the sense of traces at the top of $\Omega_0$) for all $q\in \mathbb{Z}^d$ .
\end{defin}

\begin{rem}
In the case when $s=m\in \N_0$, $u\in H^m_{\ul}(\R^d)$ if 
\begin{equation*}
\sup_{k\in \Z^d}\sum_{|\alpha|\leq m}\|\omega_{k}\partial^\alpha u\|_{L^2(\R^d)}<\infty. 
\qedhere
\end{equation*}
\end{rem}

\subsection{Solvability of the corrector equation and the two-scale expansion}

The key idea to solve~\eqref{eqn:corrector_eqns} is to replace the original Fourier multiplier $\xi \mapsto \frac{\sech(\bar h|\xi|)}{(\bar V\cdot \xi)^2-|\xi|\tanh(\bar h|\xi|)}(i\xi)$ (which is singular) by a better-behaved Fourier multiplier $\hat K_0$ which coincides with the former on $\Xi=\supp \hat b$, the only frequencies we need\footnote{{This is true for the discrete part of $\Xi$, but not quite for the non-discrete part for which we need to consider a neighborhood.}}.
We now give the appropriate notion of solution. 
\begin{defin}\label{def:notion_soln_corrector}
Let $(b,X_{T})$ satisfy the non-resonance condition of order $(0,0)$, let $K_0$ denote the convolution kernel characterized by  its Fourier transform $\hat K_0$ given in Definition~\ref{def:non-resonance}, and let $(\bar V,\bar h)\in X_T$.
For all $j\in \Z^{d}$, set 
\begin{equation}\label{eqn:psi_corrector_auxiliary}
 \psi^j_c(\cdot,\bar V, \bar h):=\bar V \cdot {K}_{0}(\cdot, \bar V, \bar h) \ast b_j, 
\end{equation}
 with $b_{j}:=\omega_{j} b$  (and $(\omega_{j})_{j\in \Z^d}$ as in Definition \ref{definition:uniformly_local}).
If $\psi_c^j$ is well-defined for all $j$, and if  $\sum_{j\in \mathbb{Z}^d}\psi_{c}^j$ is summable in $H^{1/2}_{\ul}(\R^d)$, we say that $(\psi_{c},\zeta_{c})$ given by 
\begin{equation}\label{eqn:corrector_pair}
\psi_{c}(\cdot,\bar V, \bar h):=\sum_{j\in \mathbb{Z}^d}\psi_{c}^j(\cdot,\bar V, \bar h), \quad \zeta_{c}:=-(\bar V\cdot \nabla_{y})\psi_c(\cdot,\bar V, \bar h). 
\end{equation}  
solves~\eqref{eqn:corrector_eqns}.
\end{defin}
This constitutes a good notion of solution, as the following formal computations\footnote{We shall make these computations rigorous in the proof of Theorem~\ref{thm:corrector_existence}.} suggest:
\begin{eqnarray}
\hat \psi_c(\xi)
&=& \sum_{j\in \Z^d} \bar V \cdot \hat K_{0}(\xi, \bar V, \bar h) (\hat \omega_j \ast \hat b)(\xi)\nonumber
\\ 
&=&\bar V \cdot  \hat K_{0}(\xi, \bar V, \bar h) \Big(\mathcal F (\sum_{j\in \Z^d} \omega_j) \ast \hat b\Big)(\xi)\nonumber
\\
&=&\bar V \cdot  \hat K_{0}(\xi, \bar V, \bar h)   \hat b(\xi)\nonumber
\\
&=&(\bar V \cdot i \xi) \frac{\sech(\bar h|\xi|)}{(\bar V\cdot \xi)^2-|\xi|\tanh( \bar h|\xi|)}   \hat b(\xi),\label{e.true-sol1}
\end{eqnarray}
so that $\psi_c$ satisfies \eqref{eqn:corrector_eqns} in the form of  
\begin{eqnarray}
\lefteqn{\mathcal F\Big(-\left ((\bar V\cdot \nabla_y)^2+|\nabla_{y}|\tanh{(\bar h|\nabla_y|)}\right)\psi_c\Big)(\xi)}\nonumber
\\
&=&\Big(( \bar V\cdot \xi)^2-|\xi|\tanh( \bar h|\xi|)\Big) (\bar V \cdot i \xi) \frac{\sech(\bar h|\xi|)}{(\bar V\cdot \xi)^2-|\xi|\tanh( \bar h|\xi|)}   \hat b(\xi)\nonumber
\\
&=&  (\bar V \cdot i \xi) \sech(\bar h|\xi|)  \hat b(\xi)\nonumber
\\
&=&\mathcal F \Big(( \bar V\cdot \nabla_y) \sech( \bar h|\nabla_y|)b\Big). \label{e.true-sol2}
\end{eqnarray}
The requirement that the sum converges in $H^{1/2}_{\ul}(\R^d)$ is due to the presence of the operator $\tanh(h_0|\nabla_{y}|)$ which acts on $H^{1/2}_{\ul}(\R^d)$, see also  Theorem \ref{thm:neumann}.
The following theorem establishes existence of correctors.
\begin{theor}\label{thm:corrector_existence}
Assume that $(b,X_{T})$ satisfies the non-resonance condition of order $(m,0)$ for $m\in \N_0,$  see Definition~\ref{def:non-resonance}. Then, for all $|\alpha|\leq m$ and $j\in \Z^d$, the function defined in \eqref{eqn:psi_corrector_auxiliary} satisfies for all $(\bar V,\bar h)\in X_T$
\begin{equation}\label{eqn:psi_c_j}
\partial_{y}^\alpha\psi_{c}^j(\cdot, \bar V,\bar h)=\bar V\cdot K_{\alpha}(\cdot,\bar V,\bar h) \ast b_{j} \in L^2(\R^d),
\end{equation}
where the Fourier multiplier $\hat{K}_{\alpha}$ of $K_{\alpha}$ is given in Definition~\ref{def:non-resonance}.
In addition, the pair $(\psi_{c}(\cdot, \bar V,\bar h),\zeta_{c}(\cdot, \bar V,\bar h))$ given by \eqref{eqn:corrector_pair} solves \eqref{eqn:corrector_eqns}, and belongs to $H_{\ul}^{m}(\R^d)\times H_{\ul}^{m-1}(\R^d)$.
In addition, $(y,(\bar h,\bar V)) \mapsto \psi_c(y,\bar h,\bar V) \in C^m_b(\R^d;C^m_b(X_T))$ and $(y,(\bar h,\bar V)) \mapsto \zeta_c(y,\bar h,\bar V) \in C^{m-1}_b(\R^d;C^m_b(X_T))$.
\end{theor}
Observe that only a non-resonance condition of order $(m,0)$ is required for the solvability of the corrector problem. The stronger non-resonance condition of order $(m,m)$ is required to show that the pair $(t,x) \mapsto (\psi_1,\zeta_1)(t,x,\frac x \gamma)$, with $(\psi_1,\zeta_1)$ defined as in \eqref{g00}, belongs indeed to a Sobolev space of a certain order. 
\begin{theor}\label{thm:corrector_sobolev}
Assume that $(b,X_{T})$ satisfies the non-resonance condition of order $(m,m)$ for some $m\geq 2$, see Definition \ref{def:non-resonance}, and that for some  $s>m+d/2$
\begin{equation}
(V_0,\zeta_0)\in L^\infty([0,T];H^s(\mathbb{R}^d)^{d+1})\cap {\Lip}([0,T];H^{s-1}(\mathbb{R}^d)^{d+1}).
\end{equation}
Then for all $\gamma>0$, the map $(t,x)\mapsto (\psi_{1,\gamma}, \zeta_{1,\gamma})(t,x) := (\psi_c,\zeta_c)(\frac x \gamma,V_0(t,x),h_0(t,x))$ satisfies
\begin{equation*}
\psi_{1,\gamma}\in L^\infty([0,T];H^m(\mathbb{R}^d))\cap {\Lip}([0,T];H^{m-1}(\mathbb{R}^d))
\end{equation*}
and 
\begin{equation*}
\zeta_{1,\gamma}\in L^\infty([0,T];H^{m-1}(\mathbb{R}^d))\cap {\Lip}([0,T];H^{m-2}(\mathbb{R}^d)).
\qedhere
\end{equation*}
More precisely, for all $0\le |\alpha|\le m$, $0\le |\alpha'|\le m-1$, and $0\le |\alpha''|\le m-2$,
\begin{eqnarray*}
\|\nabla^{\alpha} \psi_{1,\gamma} \|_{L^\infty([0,T];L^2\cap L^\infty(\R^d))} &\le & C \gamma^{-|\alpha|}\| (V_0,h_0)\|_{L^\infty([0,T];H^{s}(\R^d))},
\\
\|\nabla^{\alpha'} \psi_{1,\gamma} \|_{\Lip([0,T];L^2\cap L^\infty(\R^d))} &\le & C \gamma^{-|\alpha'|}\| (V_0,h_0)\|_{\Lip([0,T];H^{s-1}(\R^d))},
\\
\|\nabla^{\alpha'} \zeta_{1,\gamma} \|_{L^\infty([0,T];L^2\cap L^\infty(\R^d))} &\le & C \gamma^{-|\alpha'|}\| (V_0,h_0)\|_{L^\infty([0,T];H^{s}(\R^d))},
\\
\|\nabla^{\alpha''} \zeta_{1,\gamma} \|_{\Lip([0,T];L^2\cap L^\infty(\R^d))} &\le & C \gamma^{-|\alpha''|}\| (V_0,h_0)\|_{\Lip([0,T];H^{s-1}(\R^d))}.
\end{eqnarray*}
\end{theor}

\subsection{Proofs of Theorems~\ref{thm:corrector_existence} and~\ref{thm:corrector_sobolev}}
In order for \eqref{eqn:psi_c_j} to make sense, we need integrability conditions on the convolution kernel.
These are gathered in the upcoming proposition.
\begin{prop}\label{prop:8_bound}
Suppose that $(b,X_{T})$ satisfies the non-resonance condition of order $(m,n)\in \N^2$. 
For all $|\alpha|\leq m$ and $1\le |\beta|\leq n$, we have for all  $y\in \R^d$ and $(\bar V,\bar h) \in X_T$,
\begin{eqnarray}\label{eqn:8_bound}
|K_{\alpha}(y,\bar V,\bar h)|&\le& C
(1+|y|)^{-(d+1)} |\bar V| \|  \hat{K}_{\alpha}\|_{L^\infty_{\bar{V},\bar{h}}W^{d+1,1}_{\xi}},
\label{eqn:8_bound1}
\\
|\partial_{\bar{V},\bar{h}}^\beta K_{\alpha}(y,\bar V,\bar h)|&\le& C
(1+|y|)^{-(d+1)} \|\partial_{\bar{V},\bar{h}}^\beta \hat{K}_{\alpha}\|_{L^\infty_{\bar{V},\bar{h}}W^{d+1,1}_{\xi}}.
\label{eqn:8_bound2}
\qedhere
\end{eqnarray} 
\end{prop}
This proposition is a straightforward consequence of \eqref{eqn:prop_8_resonance} using elementary Fourier analysis: For all functions $g \in \mathcal S(\R^d)$ and all $n\in \N_0$,
\[
\sup_{x\in \R^d} |x|^n |g(x)| \,\le \, C_{d,n} \|\hat g\|_{W^{n,1}(\R^d)}.
\]
We now prove Theorems~\ref{thm:corrector_existence} and~\ref{thm:corrector_sobolev}.

\subsubsection{Proof of Theorem~\ref{thm:corrector_existence}}

Fix $(\bar V,\bar h)\in X_T$. By Proposition~\ref{prop:8_bound}, $K_{\alpha}(\cdot,\bar V,\bar h)\in L^1(\R^d)$ for all $|\alpha|\leq m$. Since $b_j=\omega_j b\in L^\infty(\R^d)$ is compactly supported, the convolutions
\[
\bar V\cdot K_{\alpha}(\cdot,\bar V,\bar h)\ast b_j
\]
are well-defined bounded continuous functions. In particular, \eqref{eqn:psi_corrector_auxiliary} defines a bounded continuous function $\psi_c^j(\cdot,\bar V,\bar h)$. Moreover, for every $|\alpha|\leq m$, we have in $\mathcal S'(\R^d)$
\[
\partial_y^\alpha \psi_c^j
=\partial_y^\alpha \big(\bar V\cdot K_0\ast b_j\big)
=\bar V\cdot \big(\partial_y^\alpha K_0\big)\ast b_j
=\bar V\cdot K_\alpha\ast b_j,
\]
because $\widehat{\partial_y^\alpha K_0}=(i\xi)^\alpha\hat K_0=\hat K_\alpha$. Since the right-hand side is a bounded continuous function, this proves \eqref{eqn:psi_c_j}; in particular $\psi_c^j(\cdot,\bar V,\bar h)\in C^m(\R^d)$.

We now establish decay properties of $\partial^\alpha \psi_c^j(\cdot,\bar V,\bar h)$. By Proposition~\ref{prop:8_bound},
\[
|K_\alpha(x,\bar V,\bar h)|\le C(1+|x|)^{-(d+1)},
\]
and since $b_j$ is supported in $j+[-1,1]^d$, we obtain for all $x\in\R^d$
\[
\begin{split}
|\partial^\alpha \psi_c^j(x,\bar V,\bar h)|
&\le |\bar V|\,\|b\|_{L^\infty(\R^d)}\int_{j+[-1,1]^d}|K_\alpha(x-z,\bar V,\bar h)|\,dz
\\
&\le C|\bar V|(1+|x-j|)^{-(d+1)}.
\end{split}
\]
This entails by summation
\[
\sup_{k\in \Z^d}  \sum_{j\in \Z^d} \|\omega_k \partial^\alpha \psi_c^j(\cdot,\bar V,\bar h)\|_{L^2(\R^d)} \,\le 
\, C |\bar V| \sup_{k\in \Z^d} \sum_{j\in \Z^d} (1+|k-j|)^{-(d+1)} \le C  |\bar V|,
\]
which shows that $\psi_c(\cdot,\bar V,\bar h) \in H^m_\ul (\R^d)$. Since $\zeta_c=-(\bar V\cdot\nabla_y)\psi_c$, we also have $\zeta_c(\cdot,\bar V,\bar h)\in H^{m-1}_\ul(\R^d)$.
Likewise, this entails that $(y,(\bar h,\bar V)) \mapsto \psi_c(y,\bar h,\bar V) \in C^m_b(\R^d;C^m_b(X_T))$ and that $(y,(\bar h,\bar V)) \mapsto \zeta_c(y,\bar h,\bar V) \in C^{m-1}_b(\R^d;C^m_b(X_T))$.

It remains to show that $(\psi_c,\zeta_c)(\cdot,\bar V,\bar h)$ solves the corrector equation~\eqref{eqn:corrector_eqns}. Since $\sum_{j\in\Z^d}\omega_j=1$ pointwise and only finitely many terms are nonzero at each point, we have $\sum_{j\in\Z^d}b_j=b$ pointwise. In addition,
\[
\sum_{j\in\Z^d}|b_j(z)|\le C|b(z)|\le C\|b\|_{L^\infty(\R^d)}
\]
for all $z\in\R^d$, where $C$ only depends on the overlap of the partition of unity. Hence, by Fubini's theorem,
\[
\begin{split}
\psi_c(x,\bar V,\bar h)
&=\sum_{j\in\Z^d}\bar V\cdot\int_{\R^d}K_0(x-z,\bar V,\bar h)b_j(z)\,dz
\\
&=\bar V\cdot\int_{\R^d}K_0(x-z,\bar V,\bar h)b(z)\,dz
=(\bar V\cdot K_0(\cdot,\bar V,\bar h))\ast b(x).
\end{split}
\]
In particular, $\psi_c$ is a bounded continuous function, hence a tempered distribution.

Let $\varphi\in\mathcal S(\R^d)$ and set
\[
\eta_\varphi(\xi):=\bar V\cdot\hat K_0(\xi,\bar V,\bar h)\,\varphi(\xi).
\]
By Proposition~\ref{prop:8_bound} and Leibniz' rule, $\eta_\varphi\in W^{d+1,1}(\R^d)$. Therefore $\widehat{\eta_\varphi}\in L^1(\R^d)$, and using $\widehat{\hat K_0}=(2\pi)^d\widetilde K_0$ with $\widetilde K_0(x,\bar V,\bar h):=K_0(-x,\bar V,\bar h)$, we obtain
\[
\widehat{\eta_\varphi}=\bar V\cdot\widetilde K_0(\cdot,\bar V,\bar h)\ast\hat\varphi.
\]
Hence,
\[
\begin{split}
\langle\widehat{\psi_c},\varphi\rangle_{\mathcal S',\mathcal S}
&=\langle\psi_c,\hat\varphi\rangle_{\mathcal S',\mathcal S}
=\int_{\R^d}\psi_c(x,\bar V,\bar h)\hat\varphi(x)\,dx
\\
&=\int_{\R^d}b(y)\big(\bar V\cdot\widetilde K_0(\cdot,\bar V,\bar h)\ast\hat\varphi\big)(y)\,dy
\\
&=\langle\hat b,\eta_\varphi\rangle_{\mathcal S',\mathcal S}.
\end{split}
\]
This proves that, in $\mathcal S'(\R^d)$, $\widehat{\psi_c}=\bar V\cdot\hat K_0(\cdot,\bar V,\bar h)\,\hat b$.

Set
\[
p(\xi):=(\bar V\cdot\xi)^2-|\xi|\tanh(\bar h|\xi|),
\qquad
q(\xi):=(\bar V\cdot i\xi)\sech(\bar h|\xi|),
\]
and
\[
\lambda(\xi):=p(\xi)\,\bar V\cdot\hat K_0(\xi,\bar V,\bar h)-q(\xi).
\]
Let $\chi\in C^\infty(\R^d)$ satisfy $\chi\equiv 1$ on a neighbourhood of $\Xi_2$ and $\supp\chi\subset \Xi_3$, and decompose
\[
\hat b=\hat b_2+\hat b_1,
\qquad
\hat b_2:=\chi\hat b,
\qquad
\hat b_1:=(1-\chi)\hat b.
\]
Since $\supp\hat b_2\subset \Xi_3$ and, by Definition~\ref{def:non-resonance},
\[
\hat K_0(\xi,\bar V,\bar h)=\frac{\sech(\bar h|\xi|)}{(\bar V\cdot \xi)^2-|\xi|\tanh(\bar h|\xi|)}(i\xi)
\qquad\text{for all }\xi\in\Xi_3,
\]
we have $\chi\lambda\equiv 0$, hence $\lambda\hat b_2=0$ in $\mathcal S'(\R^d)$. On the other hand, $\hat b_1$ has discrete support contained in $\Xi_1$, and since $b\in L^\infty(\R^d)$ the distribution $\hat b_1$ is of order $0$. By the structure theorem for order-$0$ distributions with discrete support, there exist coefficients $(a_\xi)_{\xi\in\supp\hat b_1}$ such that
\[
\hat b_1=\sum_{\xi\in\supp\hat b_1}a_\xi\delta_\xi
\]
in $\mathcal S'(\R^d)$ (the sum is locally finite). Therefore, for every $\varphi\in\mathcal S(\R^d)$,
\[
\langle \lambda\hat b_1,\varphi\rangle
=\sum_{\xi\in\supp\hat b_1}a_\xi\lambda(\xi)\varphi(\xi)=0,
\]
because $\lambda(\xi)=0$ for every $\xi\in\Xi_1$ by Definition~\ref{def:non-resonance}. Consequently,
\[
\lambda\hat b=0\qquad\text{in }\mathcal S'(\R^d).
\]
Equivalently,
\[
p(\xi)\,\widehat{\psi_c}=q(\xi)\,\hat b
\qquad\text{in }\mathcal S'(\R^d).
\]
Taking inverse Fourier transforms yields
\[
-\left ((\bar V\cdot \nabla_y)^2+|\nabla_{y}|\tanh{(\bar h|\nabla_y|)}\right)\psi_c=(\bar V\cdot \nabla_y) \sech(\bar h|\nabla_y|)b
\qquad\text{in }\mathcal S'(\R^d).
\]
Together with $\zeta_c=-(\bar V\cdot\nabla_y)\psi_c$, this proves that $(\psi_c,\zeta_c)(\cdot,\bar V,\bar h)$ solves~\eqref{eqn:corrector_eqns}.

\subsubsection{Proof of Theorem~\ref{thm:corrector_sobolev}}
We start with $\psi_{1,\gamma}(t,x):=\psi_c(\frac x\gamma,h_0(t,x),V_0(t,x))$.
On the one hand, by the assumption $(V_0,h_{0})\in L^\infty([0,T];H^s(\mathbb{R}^d,X_T))\cap {\Lip}([0,T];H^{s-1}(\mathbb{R}^d,X_T))$ for some $s>m+d/2$, we have  by the Sobolev embedding that 
\[
(V_0,h_{0})\in L^\infty([0,T];C_b^m(\mathbb{R}^d,X_T))\cap {\Lip}([0,T];C_b^{m-1}(\mathbb{R}^d,X_T)).
\]
On the other hand, by Theorem~\ref{thm:corrector_existence},
$(y,(\bar h,\bar V)) \mapsto \psi_c(y,\bar h,\bar V) \in C^m_b(\R^d;C^m_b(X_T))$.
Hence, by composition, $\psi_{1,\gamma} \in  L^\infty([0,T];C_b^m(X_T))\cap {\Lip}([0,T];C_b^{m-1}(X_T))$, and the control of these norms follow from the chain rule.
It remains to establish the square-integrability at infinity of the space-derivatives of $\psi_{1,\gamma}$.

By the proof of Theorem~\ref{thm:corrector_existence}, we have for all $j,k\in \Z^d$
\[
\sup_{k+[-1,1]^d}\sup_{\bar h} |\partial^\alpha_y \psi_c^j(\cdot,\bar V,\bar h)| \,\le\, C |\bar V|(1+|k-j|)^{-(d+1)}.
\]
Hence, 
\[
\sup_{y \in \R^d} \sup_{\bar h} |\partial^\alpha_y \psi_c(y,\bar V,\bar h)|\le\sup_{k\in \Z^d} \sum_{j\in \Z^d}\sup_{k+[-1,1]^d} \sup_{\bar h}|\partial^\alpha \psi_c^j(\cdot,\bar V,\bar h)|
\,\le \, C|\bar V|,
\]
and therefore
\begin{eqnarray*}
\int_{\R^d} |\partial^\alpha_y \psi_c(\tfrac{x}\gamma, V_0(t,x), h_0(t,x))|^2dx &\le&\int_{\R^d} \sup_{y} \sup_{\bar h} |\partial^\alpha_y \psi_c^j(y, V_0(t,x), \bar h)|^2dx\\
&\le& C \int_{\R^d} |V_0(t,x)|^2 dx.
\end{eqnarray*}
For $\alpha=0$, this shows that $\psi_{1,\gamma} \in L^\infty([0,T];L^2(\R^d))$
and 
\[ \|\psi_{1,\gamma}\|_{L^\infty([0,T];L^2(\R^d))} \le C \| V_0\|_{L^\infty([0,T];L^2(\R^d))} .\]
By the chain rule and the full strength of Proposition~\ref{prop:8_bound}, the same argument entails that $\psi_{1,\gamma} \in L^\infty([0,T];H^m(\R^d)) \cap \Lip ([0,T];H^{m-1}(\R^d))$.
{In addition, for each derivative with respect to $x$, there is a scaling factor $\gamma^{-1}$, and we have more generally}
for all $0\le |\alpha|\le m$
\[
\|\nabla^\alpha \psi_{1,\gamma} \|_{L^\infty([0,T];L^2(\R^d))} \,\le \, C \gamma^{-|\alpha|}\| (V_0,h_0)\|_{L^\infty([0,T];H^{|\alpha|}(\R^d))}.
\]
Likewise, for all $0\le |\alpha|\le m-1$,
\[
\|\nabla^\alpha \psi_{1,\gamma} \|_{\Lip([0,T];L^2(\R^d))} \,\le \, C \gamma^{-|\alpha|}\| (V_0,h_0)\|_{\Lip([0,T];H^{|\alpha|}(\R^d))}.
\]
Since $\zeta_{1,\gamma}(t,x)=-V_0\cdot (\nabla_{y}\psi_{c})|_{y=\frac x\gamma,\bar h=h_0(t,x),\bar V=V_0(t,x)}$, the corresponding claim follows by the chain rule.

%%%%%%%%
%%%%%%%%

\section{The two-scale Dirichlet-to-Neumann map}

Now that we have constructed correctors, we have to characterize the Dirichlet-to-Neumann map for the elliptic problem \eqref{g.7}.

\subsection{Dirichlet-to-Neumann map}
We rewrite  \eqref{g.7}  as 
\begin{equation}\label{eqn:elliptic_uniformly_local}
\left\{
\begin{array}{rcl}    
(h_0^2\Delta_{y}+\partial_{z}^2)\Phi&=&0 , \quad \text{on $\Omega_0$},\\
\Phi|_{z=0}&=&\Theta,\\
\frac{1}{h_0}\partial_{z}\Phi|_{z=-1} &=&\Gamma,
\end{array}
\right.
\end{equation}
where $\Theta=\psi_1$ and $\Gamma=\nabla_{y} b  \cdot \nabla \psi_0$. {Note that the slow variables $(t,x)$ are frozen, so $h_0$ is fixed. By assumption, $h_0$ ranges in a compact subset of $(0,\infty)$. In what follows, multiplicative constants may depend on $h_0$ but remain uniformly controlled on that compact set.} We will refer to $(\Theta,\Gamma)$ as the boundary data of~\eqref{eqn:elliptic_uniformly_local}.  

\medskip

We define the Neumann map   as 
\begin{equation}\label{eqn:lambda}
\calN:(\Theta,\Gamma)\mapsto \frac{1}{h_0}\partial_{z}\Phi|_{z=0},
\end{equation}
where $\Phi$ is the solution of \eqref{eqn:elliptic_uniformly_local} with boundary data $(\Theta,\Gamma)$,  see Definition \ref{def:elliptic} for the precise notion of solution. The following result expresses the Neumann map as
a pseudo-differential operator.
\begin{theor0}\label{thm:neumann}
Assume that $(\Theta,\Gamma)\in H^{1/2}_{\ul}(\R^d)\times H^{-1/2}_{\ul}(\R^d)$. Then, \eqref{eqn:elliptic_uniformly_local}  has a unique solution. Moreover,  the Neumann map $\calN:H^{1/2}_{\ul}(\R^d)\times H^{-1/2}_{\ul}(\R^d)\rightarrow H^{-1/2}_{\ul}(\R^d)$, defined in \eqref{eqn:lambda}, is continuous and can be alternatively expressed by
\begin{equation}\label{eqn:neumann_mapping}
\calN(\Theta,\Gamma)=|\nabla_{y}|\tanh(h_0|\nabla_{y}|)\Theta+\sech(h_0|\nabla_{y}|)\Gamma. 
\end{equation}
\end{theor0}
\begin{rem}
A version of Theorem \ref{thm:neumann} can be found in \cite[Section 3.3.1]{MR2901196} for  $\Theta$ and $\Gamma$  periodic. In this particular case, Fourier series can be used to obtain explicit expressions for the solution and for the Neumann map. Our approach is inspired by~\cite[Section 3]{MR3465379}. This result of Alazard, Burq and Zuily considers more general domains. Our situation is more specific and our main focus is the explicit identification  \eqref{eqn:neumann_mapping} of the Neumann map. 
\end{rem}

\begin{rem}
The precise meaning of $|\nabla_{y}|\tanh(h_0|\nabla_{y}|)\Theta+\sech(h_0|\nabla_{y}|)\Gamma$ in  \eqref{eqn:neumann_mapping} is 
\begin{equation}\label{eqn:sum_tanh_sech}
	|\nabla_{y}|\tanh(h_0|\nabla_{y}|)\Theta+\sech(h_0|\nabla_{y}|)\Gamma=\sum_{q\in \Z^d} 	|\nabla_{y}|\tanh(h_0|\nabla_{y}|)(\omega_q\Theta)+\sech(h_0|\nabla_{y}|)(\omega_q\Gamma),
\end{equation}
where $(\omega_{q})_{q\in \Z^d}$ is given  in Definition~\ref{definition:uniformly_local} and the summand is understood  for all $q\in \mathbb{Z}^d$ as
\begin{equation}
\begin{split}
	|\nabla_y|\tanh(h_0|\nabla_{y}|)(\omega_q\Theta)&=\mathcal{F}^{-1}(|\cdot|\tanh(h_0|\cdot|))*(\omega_{q}\Theta), \\	\sech(h_0|\nabla_{y}|)(\omega_q\Gamma)&=\mathcal{F}^{-1}(\sech(h_0|\cdot|))*(\omega_{q}\Gamma),
\end{split}
\end{equation}
and well-defined by our assumption $(\Theta,\Gamma)\in H^{1/2}_{\ul}(\R^d)\times H^{-1/2}_{\ul}(\R^d)$.
\end{rem}

\subsection{Proof of Theorem \ref{thm:neumann}}
Let us first outline the strategy. Since $(\Theta,\Gamma)$ are not necessarily integrable, we need an appropriate notion of solution, which we give in Definition~\ref{def:elliptic} below. To do so, we localize the boundary data and study the following elliptic problem for all $q\in \Z^d$ 
\begin{equation}\label{eqn:elliptic_uniformly_localized}
\left\{
\begin{array}{rcl}    
(h_0^2\Delta_{y}+\partial_{z}^2)\Phi_{q}&=&0 , \quad \text{on $\Omega_0$},\\
\Phi_{q}|_{z=0}&=&\Theta_{q},\\
\frac{1}{h_0}\partial_{z}\Phi_{q}|_{z=-1} &=&\Gamma_{q},
\end{array}
\right.
\end{equation}
where $\Theta_{q}:=\omega_{q}\Theta$, $\Gamma_{q}:=\omega_{q}\Gamma$, and $\omega_{q}$  as in Definition~\ref{definition:uniformly_local}.
To prove Theorem~\ref{thm:neumann}, we  need to show that $\Phi:=\sum_{q\in \Z^d}\Phi_{q}$ indeed solves
\eqref{eqn:elliptic_uniformly_local}. From this, the characterization~\eqref{eqn:neumann_mapping} follows from the identity 
\begin{equation}\label{eqn:DTN_characterization_2}
\frac{1}{h_0}\partial_{z} \Phi_{q}|_{z=0}=|\nabla_{y}|\tanh (h_0|\nabla_{y}|)\Theta_{q} +\sech(h_0|\nabla_{y} |)\Gamma_{q}.
\end{equation}
The key estimate\footnote{Recall the notation in \eqref{e.not-nab-mu}.} to show that $\Phi$ belongs to a uniformly local Sobolev space is the following adaptation of~\cite[Lemma 3.6]{MR3465379}: there exists some $\delta>0$ such that for all $q\in \mathbb{Z}^d$
\begin{equation}\label{eqn:ulss_est}
\|e^{\delta\langle \cdot-q \rangle }\nabla^{h_0^2}\Phi_{q}\|_{L^2(\Omega_0)}\leq C ( \|\Theta_{q}\|_{H^{1/2}(\mathbb{R}^d)}+\|\Gamma_{q}\|_{H^{-1/2}(\mathbb{R}^d)}),
\end{equation}
where $C>0$ is a constant independent\footnote{This constant depends on $h_0$.} of $q$, which we establish  in Proposition~\ref{prop:bounds_variational_soln}. 

Before we define the notion of solution, we first need to extend $\Theta\in H^{1/2}_{\ul}(\R^d)$ into an element $\Theta^\dagger\in H^1_{\ul}(\Omega_0)$. The existence of such an extension follows from \cite[Lemma 3.5]{MR3465379} which we briefly outline here. For all $q\in \mathbb{Z}^d$, one can extend $\Theta_{q}=\omega_q \Theta$ on $\Omega_0$  as $\Theta_{q}^\dagger$ in such a way that $\Theta_{q}^\dagger\in H^1(\Omega_0)$,
\begin{equation}\label{eqn:support_extension}
\supp(\Theta_{q}^\dagger)\subset \{(x,z)\in \Omega_0:|x-q|\leq 2,z\in (-1/2,0]\},
\end{equation}
and, for some $C>0$ that is independent of $q$,
\begin{equation}\label{eqn:trace_ineq}
\|\Theta_{q}^\dagger\|_{H^1(\Omega_0)}\leq C\|\Theta_{q}\|_{H^{1/2}(\mathbb{R}^d)}.
\end{equation}
We then define $\Theta^\dagger:=\sum_{q\in \mathbb{Z}^d}\Theta_{q}^\dagger$, which, by the above, indeed satisfies $\Theta^\dagger\in H^1_{\ul}(\Omega_0)$ and $\|\Theta^\dagger\|_{H^1_{\ul}(\Omega_0)}\leq  C \|\Theta\|_{H^{1/2}_{\ul}(\R^d)}$.

\medskip

We are now in position to define a notion of solution for~\eqref{eqn:elliptic_uniformly_local}. 
\begin{defin}\label{def:elliptic}
 If there exists $\tilde{\Phi}\in H_{0,\ul}^{1}(\Omega_0)$ such that\footnote{Recall that $H^1_0(\Omega_0):=\{\psi \in H^1(\Omega_0)\,|\, \psi_{|z=0}\equiv 0 \text{ in the sense of traces}\}$.} for all  $\psi \in H^1_0(\Omega_0)$
\begin{equation}\label{eqn:definition_soln}
\int_{\Omega_0}\nabla^{h_0^2}  \psi\cdot\nabla^{h_0^2} \tilde{\Phi}  +\int_{\Omega_0}\nabla^{h_0^2}  \psi\cdot\nabla^{h_0^2} \Theta^\dagger    +h_0 \langle \Gamma, \psi(\cdot,-1) \rangle_{H^{-1/2}(\R^d),H^{1/2}(\R^d)} =0,
\end{equation}
then we say that $\Phi:=\tilde{\Phi}+\Theta^\dagger$ is a solution to \eqref{eqn:elliptic_uniformly_local}.
\end{defin}
In what follows, we shall abusively write 
\[
\langle f, g(\cdot,-1) \rangle_{H^{-1/2}(\R^d),H^{1/2}(\R^d)} =\int_{\mathbb{R}^d} f g|_{z=-1}.
\]
The proof of uniqueness of \eqref{def:elliptic} is the same as in \cite[Proposition 3.3]{MR3465379}, and we only argue in favor of existence of a solution, and establish the explicit formula of the Neumann map. 
Consider the standard elliptic problem \eqref{eqn:elliptic_uniformly_localized} with localized boundary data. 
By Lax-Milgram's theorem, there exists a unique $\tilde{\Phi}_{q}\in H_{0}^{1}(\Omega_0)$
such that for all $\psi \in H_{0}^{1}(\Omega_0)$
\begin{equation}\label{eqn:elliptic_variational_weak}
\int_{\Omega_0}\nabla^{h_0^2}\psi\cdot \nabla^{h_0^2} \tilde{\Phi}_q    =-\int_{\Omega_0}\nabla^{h_0^2}\psi\cdot \nabla^{h_0^2} \Theta_q^\dagger    -h_0\int_{\mathbb{R}^d}\Gamma_q \psi |_{z=-1} .
\end{equation}
Hence, $\Phi_{q}=\Theta_{q}^{\dagger}+\tilde{\Phi}_{q}$ is the unique weak solution of~\eqref{eqn:elliptic_uniformly_localized}. The following proposition, which we prove in Section~\ref{subsec:proof_prop_variational_bdds}, then implies Theorem~\ref{thm:neumann}. 
\begin{prop}\label{prop:bounds_variational_soln}
The unique weak solution $\Phi_{q}\in H^1_0(\Omega_0)$ of \eqref{eqn:elliptic_uniformly_localized} satisfies
for all $z\in [-1,0]$
\begin{equation}\label{eqn:dtn_characterization_1}
\Phi_{q}(\cdot,z)=\frac{\cosh(h_0(z+1)|\nabla_{y}|)}{\cosh(h_0|\nabla_{y}|)}
\Theta_{q} + \frac{\sinh(h_0z |\nabla_{y}|)}{|\nabla_{y}|\cosh(h_0|\nabla_{y}|)}
\Gamma_{q},
\end{equation}
as well as 
\begin{equation}\label{eqn:dtn_characterization_2}
\frac{1}{
h_0}\partial_{z}\Phi_{q}|_{z=0}=|\nabla_{y}|\tanh(h_0|\nabla_{y}|)\Theta_{q}+\sech(h_0|\nabla_{y}|)\Gamma_{q}. 
\end{equation}
In addition, there exist some $\delta>0$ and some constant $C>0$ (both depending on $h_0$ but independent of $q$) such that  	
\begin{equation}\label{eqn:cauchy_seq}
\|e^{\delta\langle \cdot-q \rangle }\nabla^{h_0^2} \Phi_{q}\|_{L^2(\Omega_0)}\leq  C(\|\Theta_{q}\|_{H^{1/2}(\mathbb{R}^d)}+\|\Gamma_{q}\|_{H^{-1/2}(\mathbb{R}^d)}).
\end{equation}
\end{prop}
We split the rest of this section in two paragraphs. 
We first prove Theorem~\ref{thm:neumann} based on Proposition~\ref{prop:bounds_variational_soln},
and then turn to the proof of Proposition~\ref{prop:bounds_variational_soln}.

\subsubsection{Proof of Theorem \ref{thm:neumann}}
By density it is enough to assume that $\Gamma_q,\Theta_q \in \mathcal S(\R^d)$. We proceed in three steps.

\medskip

\step1 Well-posedness of \eqref{eqn:elliptic_uniformly_local}.\\
Uniqueness is as for~\cite[Proposition 3.3]{MR3465379}, and we only prove existence. 
By \eqref{eqn:cauchy_seq} in Proposition~\ref{prop:bounds_variational_soln}, we have for all $k,q\in \mathbb{Z}^d$ 
\begin{equation*}
\begin{split}
 \|\omega_{k}\nabla^{h_0^2} {\Phi}_{q}  \|_{L^2(\Omega_0)}&=  \|\omega_{k}e^{-\delta \langle \cdot-q\rangle }e^{\delta \langle \cdot-q\rangle }\nabla^{h_0^2} {\Phi}_{q}  \|_{L^2(\Omega_0)},\\
&\leq Ce^{-\delta \langle k-q\rangle }\Big(\|\Theta_{q}\|_{H^{1/2}(\mathbb{R}^d)}+
\|\Gamma_{q}\|_{H^{-1/2}(\mathbb{R}^d)}\Big),
\end{split}
\end{equation*}
for some constant $C>0$ that is independent of $k,q\in \Z^d$. 
This implies that $\Phi:=\sum_{q\in \Z^d} \Phi_q$ satisfies
\begin{equation*}
\begin{split}
 \|\omega_{k}\nabla^{h_0^2}{\Phi}  \|_{L^2(\Omega_0)}
&\leq \sum_{q\in \mathbb{Z}^d} \|\omega_{k}\nabla^{h_0^2} {\Phi}_{q}  \|_{L^2(\Omega_0)},\\
&\leq C\sum_{q\in \mathbb{Z}^d}e^{-\delta \langle k-q\rangle }\Big(\|\Theta_{q}\|_{H^{1/2}(\mathbb{R}^d)}+
\|\Gamma_{q}\|_{H^{-1/2}(\mathbb{R}^d)}\Big)\\
&\leq C\Big(\|\Theta_{q}\|_{H^{1/2}_\ul(\mathbb{R}^d)}+
\|\Gamma_{q}\|_{H^{-1/2}_\ul(\mathbb{R}^d)}\Big)\sum_{q\in \mathbb{Z}^d}e^{-\delta \langle k-q\rangle }.
\end{split}
\end{equation*}
Hence, $\nabla\Phi \in L^2_\ul(\Omega_0)$ and by (local) Poincar\'e's inequality, $\Phi \in H^1_\ul(\Omega_0)$, and 
\begin{equation*}
 \| \Phi \|_{H^1_{\ul}(\Omega_0)}\leq C \left(\|\Theta\|_{H^{1/2}_{\ul}(\mathbb{R}^d)}+ \|\Gamma\|_{H^{-1/2}_{\ul}(\mathbb{R}^d)}\right ).
\end{equation*}
We now only have to check that our choice of $\Phi$ is a solution in the sense of Definition \ref{def:elliptic}. 
Recall that $\tilde{\Phi}=\Phi-\Theta^\dagger=\sum_{q\in \mathbb{Z}^d}\tilde{\Phi}_{q}$, where $\tilde{\Phi}_q$ satisfies~\eqref{eqn:elliptic_variational_weak}. Hence, for all 
$\psi \in H^{1}_0(\Omega_0)$ such that $\psi|_{|y|>K} \equiv 0$ for some $K>0$, we have  
\begin{equation*}
\begin{split}
\int_{\Omega_0}\nabla^{h_0} \tilde{\Phi} \cdot \nabla^{h_0^2} \psi&=
\sum_{q\in \Z^d}\int_{\Omega_0}\nabla^{h_0^2} \tilde{\Phi}_q \cdot \nabla^{h_0^2}  \psi \\
&=-\sum_{q\in \Z^d}\int_{\Omega_0}\nabla^{h_0^2} \Theta_q^\dagger \cdot\nabla^{h_0^2} \psi +h_0\int_{\mathbb{R}^d}\Gamma_q  \psi|_{z=-1}\\
&=-\int_{\Omega_0}\nabla^{h_0^2} \Theta^\dagger \cdot \nabla^{h_0^2}  \psi -h_0\int_{\mathbb{R}^d}\Gamma  \psi |_{z=-1} ,
\end{split}
\end{equation*}
where the sum on $q\in \Z^d$ is on a finite set (by the support condition on $\psi$). By letting $K\uparrow +\infty$ and density argument, this shows that $\Phi$ solves \eqref{eqn:elliptic_uniformly_local}. 

\medskip

\step2 Continuity of the Neumann map.\\
To prove this continuity of the Neumann map, we fix $k\in \Z^d$ and show that 
\begin{equation*}
\left \| \omega_{k}\partial_{z}\Phi|_{z=0} \right \|_{H^{-1/2}(\mathbb{R}^d)}\leq C\left (\|\Theta\|_{H^{1/2}_{\ul}(\mathbb{R}^d)}+\|\Gamma\|_{H^{-1/2}_{\ul}(\mathbb{R}^d)} \right ),
\end{equation*}
where $C>0$ is independent of $k\in \Z^d$.
By a trace estimate and the triangle inequality,
\begin{equation*}
 \| \omega_{k}\partial_{z}\Phi|_{z=0}  \|_{H^{-1/2}(\mathbb{R}^d)}
\le C   \| \omega_{k}\partial_{z}\Phi  \|_{L^2(\Omega_0)} \le C \sum_{q \in \Z^d} 
\| \omega_{k}\nabla^{h_0^2}\Phi_q  \|_{L^2(\Omega_0)} .
\end{equation*}
Moreover, by Proposition~\ref{prop:bounds_variational_soln}  and the support condition on $\omega_{k}$, we have  
\begin{equation*}
\begin{split}
\sum_{q\in \mathbb{Z}^d}  \| \omega_{k}\nabla^{h_0^2}\Phi_{q}   \|_{L^2(\Omega_0)}&= \sum_{q\in \mathbb{Z}^d}  \| \omega_{k}e^{-\delta \langle \cdot-q\rangle}e^{\delta \langle \cdot-q\rangle}\nabla^{h_0^2}\Phi_{q}   \|_{L^2(\Omega_0)},\\
&\leq C\sum_{q\in \mathbb{Z}^d}e^{-\delta\langle k-q\rangle }  \|e^{\delta \langle \cdot-q\rangle}\nabla^{h_0^2}\Phi_{q}   \|_{L^2(\Omega_0)}, \\
&\leq C\sum_{q\in \mathbb{Z}^d}e^{-\delta\langle k-q\rangle }(\|\Theta_{q}\|_{H^{1/2}(\mathbb{R}^d)}+\|\Gamma_{q}\|_{H^{-1/2}}),\\
&\leq C\left (\|\Theta\|_{H^{1/2}_{\ul}(\mathbb{R}^d)}+\|\Gamma\|_{H^{-1/2}_{\ul}(\mathbb{R}^d)}\right )\sum_{q\in \mathbb{Z}^d}e^{-\delta\langle k-q\rangle}.
\end{split}
\end{equation*}
Combining these inequalities yields 
\begin{equation*}
\left \| \partial_{z}\Phi|_{z=0} \right \|_{H_{\ul}^{-1/2}(\mathbb{R}^d)}
=\sup_{k\in \Z^d}\left \| \omega_{k}\partial_{z}\Phi|_{z=0} \right \|_{H^{-1/2}(\mathbb{R}^d)}\leq C\left (\|\Theta\|_{H^{1/2}_{\ul}(\mathbb{R}^d)}+\|\Gamma\|_{H^{-1/2}_{\ul}(\mathbb{R}^d)} \right ).
\end{equation*}
Hence, we have that the Neumann map
\begin{equation*}
\begin{array}{rrcl}
\calN:&H^{1/2}_{\ul}(\mathbb{R}^d)\times H^{-1/2}_{\ul}(\mathbb{R}^d)&\rightarrow &H^{-1/2}_{\ul}(\mathbb{R}^d),\\
&(\Theta,\Gamma)&\mapsto& \frac{1}{h_0}\partial_{z}\Phi|_{z=0},
\end{array}
\end{equation*}
is continuous. 

\medskip

\step3 Formula for the Neumann map.\\
By Proposition \ref{prop:bounds_variational_soln} and continuity of the Neumann map, we have  
\begin{equation*}
\begin{split}
\calN(\Theta,\Gamma)&= \sum_{q\in \Z^d}\calN(\Theta_{q},\Gamma_{q}),\\
&=\sum_{q\in \Z^d}\frac{1}{h_0}\partial_{z}\Phi_{q}|_{z=0},\\
&=\sum_{q\in \Z^d}|\nabla_{y}|\tanh(h_0|\nabla_{y}|)\Theta_{q}+\sech(h_0|\nabla_{y}|)\Gamma_{q},\\
&=|\nabla_{y}|\tanh(h_0|\nabla_{y}|)\Theta+\sech(h_0|\nabla_{y}|)\Gamma,
\end{split}
\end{equation*}
as claimed.

\subsubsection{Proof of Proposition~\ref{prop:bounds_variational_soln}} \label{subsec:proof_prop_variational_bdds}

We split the proof into two steps. By density, we can assume that $\Theta_q,\Gamma_q \in  \mathcal{S}(\mathbb{R}^d)$. 

\medskip

\step1 Proof of~\eqref{eqn:dtn_characterization_1} and~\eqref{eqn:dtn_characterization_2}.
\\
Taking the Fourier transform of~\eqref{eqn:elliptic_uniformly_localized} with respect to the horizontal variables, we have that for all $\xi \in \R^d$ and $z\in (-1,0)$,
\begin{equation*}
\begin{cases}
&(-h_0^2|\xi|^2+\partial_{z}^2)\hat{\Phi}_{q}(\xi,z)=0 ,\\
&\hat{\Phi}_{q}(\xi,0)=\hat{\Theta}_{q}(\xi),\\
&\partial_{z} \hat{\Phi}_q(\xi,-1)=h_0\hat{\Gamma}_{q}(\xi). 
\end{cases}
\end{equation*}
This is the Sturm-Liouville problem of \cite[Proposition 3.5]{MR2901196}, which one can explicitly solve: For all $(\xi,z)\in \mathbb{R}^d\times[-1,0]$,
\begin{equation}\label{eqn:explicit_soln}
\hat{\Phi}_{q}(\xi,z)=\frac{\cosh(h_0(z+1)|\xi|)}{\cosh(h_0|\xi|)}
\hat{\Theta}_{q}(\xi) + \frac{1}{\cosh(h_0|\xi|)} \frac{\sinh(zh_0|\xi|)}{|\xi|}
\hat{\Gamma}_{q}(\xi),\end{equation}
which yields upon differentiating at $z=0$ 
\begin{equation*}
\frac{1}{h_0}\partial_{z}\hat{\Phi}_{q}(\xi,0)=|\xi|\tanh(h_0|\xi|)\hat{\Theta}_{q}(\xi)+ \sech(h_0|\xi|)\hat{\Gamma}_{q}(\xi).
\end{equation*}
In \eqref{eqn:explicit_soln}, we identified $\xi\mapsto \frac{\sinh(zh_0|\xi|)}{|\xi|}$ 
 with its smooth extension. By Fourier inversion, this yields the operator forms \eqref{eqn:dtn_characterization_1} and \eqref{eqn:dtn_characterization_2}. 

\medskip

\step2 Proof of~\eqref{eqn:cauchy_seq}.\\
Without loss of generality, we assume that $\delta\in (0,1)$.
Since $\Phi_q=\tilde{\Phi}_q+\Theta^\dagger_q$,  
\begin{equation}\label{ab1}
\|e^{\delta \langle \cdot-q\rangle }\nabla^{h_0^2}\Phi_q\|_{L^2(\Omega_0)}\leq \|e^{\delta \langle \cdot-q\rangle }\nabla^{h_0^2}\tilde{\Phi}_q\|_{L^2(\Omega_0)}+\|e^{\delta \langle \cdot-q\rangle }\nabla^{h_0^2}\Theta^\dagger_q\|_{L^2(\Omega_0)},
\end{equation}
and we bound each term separately. 
The control of the second term is easy.
By the support condition~\eqref{eqn:support_extension} on $\Theta_{q}^\dagger$,  we have
$e^{\delta \langle \cdot-q\rangle }\leq  e^{\sqrt{5}}$ on $\supp(\Theta_{q}^\dagger)$.
Hence, by~\eqref{eqn:trace_ineq},
\begin{equation}\label{eqn:bdd_dagger}
\|e^{\delta \langle \cdot-q\rangle }\nabla^{h_0^2}\Theta^\dagger_q\|_{L^2(\Omega_0)}\le 
C\|\nabla^{h_0^2}\Theta^\dagger_q\|_{L^2(\Omega_0)}\le  C \|\Theta_{q}\|_{H^{1/2}(\mathbb{R}^d)}.
\end{equation}
To treat the first term on the right-hand side of~\eqref{ab1}, we need two technical ingredients.
First, we recall the following weighted Poincar\'e inequality of \cite[Remark 3.2]{MR3465379}.
There exists $C>0$ such that for all measurable weights $\alpha\in C^\infty_{b}(\mathbb{R}^d,\R_+)$ and all functions $u\in H^{1}_{0}(\Omega_0)$ we have\footnote{The proof is elementary, and amounts to writing $u(x,z')=-\int_{z'}^0 \partial_z u(x,z) dz$ and using that $\alpha$ is independent of $z$.} 
\begin{equation}\label{eqn:poincare}
\int_{\Omega_0}\alpha(x)|u(x,z)|^2\,dx\,dz\leq C \int_{\Omega_0}\alpha(x)|\partial_{z}u(x,z)|^2\,dx\,dz.
\end{equation}
Second, we introduce a bounded and Lipschitz-continuous approximation $\eta_{\e, q}$ of $\langle \cdot -q\rangle$ parametrized by $0<\e<1$ and defined for all $x\in \R^d$ by
\begin{equation}
\eta_{\e, q}(x):=\frac{\langle x-q \rangle}{1+\e\langle x-q \rangle},
\end{equation}
and which satisfies
\begin{equation}\label{lem:eta}
\|\eta_{\e,q}\|_{L^\infty(\mathbb{R}^d)}\leq \e^{-1}\text{ and } \|\nabla \eta_{\e,q}\|_{L^\infty(\mathbb{R}^d)}\leq 1.
\end{equation}
In particular, $e^{\delta \eta_{\e,q}}\in W^{1,\infty}(\mathbb{R}^d)$ for $\delta>0$.  

We now proceed with the estimate of the first right-hand side term of \eqref{ab1}. 
By monotone convergence, it is enough to show that for $0<\delta<1$ small enough, we have 
for all $0<\e <1$
\begin{equation}\label{eqn:bdd_tilde}
\left (\int_{\Omega_0}e^{2\delta \eta_{\e,q}}|\nabla^{h_0^2} \tilde{\Phi}_q|^2 \right )^\frac{1}{2}\le C( \|\Theta_{q}\|_{H^{1/2}(\mathbb{R}^d)}+\left \| \Gamma_q \right \|_{H^{-1/2}(\mathbb{R}^d)}). 
\end{equation}
By construction\footnote{This is the reason why we introduced $\eta_{\e,q}$.}, with $\chi_{\delta,\e} := e^{\delta \eta_{\e,q}}$,  $\psi := \chi_{\delta,\e}^2\tilde{\Phi}_{q}$ is of class $H_{0}^{1}(\Omega_0)$, so that the integral identity~\eqref{eqn:elliptic_variational_weak} holds, and we may proceed as for Caccioppoli's estimate. We fix the weight $\alpha= \chi_{\delta,\e}^2$. The first two terms are standard.
On the one hand,
\begin{eqnarray*}
\int_{\Omega_0} \nabla^{h_0^2} \psi \cdot  \nabla^{h_0^2}\tilde{\Phi}_{q}&\ge & \int_{\Omega_0}\chi_{\delta,\e}^2 |\nabla^{h_0^2} \tilde{\Phi}_{q}|^2- C\delta^2 \int_{\Omega_0} (\chi_{\delta,\e} \tilde{\Phi}_{q})^2 |\nabla^{h_0^2}  \eta_{\e,q}|^2,
\end{eqnarray*}
which, by Poincar\'e's inequality~\eqref{eqn:poincare}, yields for $0<\delta<1$ small enough\footnote{This is the only place we use the smallness of $\delta$.}
\[
\int_{\Omega_0} \nabla^{h_0^2} \psi \cdot  \nabla^{h_0^2}\tilde{\Phi}_{q}\,\ge \, \frac34 \int_{\Omega_0}\chi_{\delta,\e}^2 |\nabla^{h_0^2} \tilde{\Phi}_{q}|^2.
\]
On the other hand, by Poincar\'e's inequality~\eqref{eqn:poincare} and by \eqref{eqn:bdd_dagger},
\begin{eqnarray*}
\lefteqn{\Big|\int_{\Omega_0} \nabla^{h_0^2} \psi \cdot  \nabla^{h_0^2}\Theta^\dagger_q\Big|}
\\
&\le&C \Big( \int_{\Omega_0}\chi_{\delta,\e}^2 |\nabla^{h_0^2}\tilde{\Phi}_{q}|^2+\delta^2 \chi_{\delta,\e}^2 \tilde{\Phi}_{q}^2 |\nabla^{h_0^2}  \eta_{\e,q}|^2 \Big)^\frac12\Big(\int_{\Omega_0} \chi_{\delta,\e}^2 | \nabla^{h_0^2}\Theta^\dagger_q|^2\Big)^\frac12
\\
&\le &C\Big( \int_{\Omega_0}\chi_{\delta,\e}^2 |\nabla^{h_0^2} \tilde{\Phi}_{q}|^2\Big)^\frac12 \|\Theta_{q}\|_{H^{1/2}(\mathbb{R}^d)}.
\end{eqnarray*}
The only subtle term is due to the Neumann boundary condition.
Let $0\le \bar \omega_q\le 1$ be a smooth function such that $\bar \omega_q |_{\supp(\omega_q)} \equiv 1$ and
\begin{equation*}
\supp(\bar \omega_q)\subset \{(x,z)\in \Omega_0:|x-q|\leq 3,z\in (-2/3,0]\}, \quad \sup |\nabla \bar \omega_q| \le C.
\end{equation*}
By the support condition on $\Gamma_q$, we have
\[
\Big|\int_{\R^d} \chi_{\delta,\e}^2 \Gamma_q \tilde{\Phi}_{q}|_{z=-1}\Big| =\Big| \int_{\R^d} \bar \omega_q \chi_{\delta,\e}^2 \Gamma_q \tilde{\Phi}_{q}|_{z=-1} \Big|\le  \|\Gamma_q\|_{H^{-1/2}(\R^d)} \|\bar \omega_q \chi_{\delta,\e}^2  \tilde{\Phi}_{q}\|_{H^{1/2}(\R^d)}.
\]
By trace theory and Poincar\'e's inequality~\eqref{eqn:poincare}, 
we can control the second factor by
\begin{eqnarray*}
\lefteqn{\|\bar \omega_q \chi_{\delta,\e}^2  \tilde{\Phi}_{q}\|_{H^{1/2}(\R^d)}}
 \\
&\le&C  
 \|\bar \omega_q \chi_{\delta,\e}^2   \tilde{\Phi}_{q}\|_{H^{1}(\Omega_0)}
\\
&\le & C \Big(  \|\bar \omega_q \chi_{\delta,\e}^2   \tilde{\Phi}_{q}\|_{L^2(\Omega_0)}
+ \|\bar \omega_q \chi_{\delta,\e}^2   \nabla \tilde{\Phi}_{q}\|_{L^2(\Omega_0)}
+\| (\chi_{\delta,\e} \nabla\bar \omega_q +2\bar \omega_q \nabla \chi_{\delta,\e}) \chi_{\delta,\e}    \tilde{\Phi}_{q}\|_{L^2(\Omega_0)}
\Big)
\\
&\le &C \|  \bar \omega_q \chi_{\delta,\e}\|_{W^{1,\infty}(\Omega_0)}
\Big(\| \chi_{\delta,\e} \tilde{\Phi}_{q}\|_{L^2(\Omega_0)}+\| \chi_{\delta,\e} \nabla \tilde{\Phi}_{q}\|_{L^2(\Omega_0)}\Big)
\\
&\le &C  \|\chi_{\delta,\e} \nabla^{h_0^2} \tilde{\Phi}_{q}\|_{L^2(\Omega_0)},
\end{eqnarray*}
where we used that $|\chi_{\delta,\e}|+|\nabla \chi_{\delta,\e}|$ is uniformly bounded with respect to $\delta$ on $\supp(\bar \omega_q)$.
All in all, we have thus proved 
\[
 \|\chi_{\delta,\e} \nabla^{h_0^2} \tilde{\Phi}_{q}\|_{L^2(\Omega_0)}^2
\le C \|\chi_{\delta,\e} \nabla^{h_0^2} \tilde{\Phi}_{q}\|_{L^2(\Omega_0)}\Big( \|\Theta_{q}\|_{H^{1/2}(\mathbb{R}^d)} + \|\Gamma_q\|_{H^{-1/2}(\R^d)} \Big) ,
\]
which yields \eqref{eqn:bdd_tilde}, and concludes the proof.

\section{Proof of the consistency result}\label{sec:proof}

\subsection{Structure of the proof}

Recall that we use the definition \eqref{g00} of $\zeta_1$ and $\psi_1$
\begin{equation*} 
\zeta_1(t,x,y)=\zeta_c(y;(h_0(t,x), V_0(t,x))), \quad \psi_1(t,x,y)=\psi_c(y;(h_0(t,x), V_0(t,x))),
\end{equation*}
and the definition~\eqref{eqn:DTN_ref} \&~\eqref{eqn:DTN:eff} of $G_{\text{res}}$.
To simplify the presentation we omit the dependence on the time variable.
As a starting point, we shall establish that the error terms in Theorem~\ref{thm:main} can be reformulated as
\begin{equation}\label{eqn:error_1}
\begin{split}
E_{1}\left (\zeta_{\mu}^{2s},\psi_{\mu}^{2s}\right )&=\sqrt{\mu}(\partial_{t}\zeta_1-G_{\text{res}}),\\
E_{2}\left (\zeta_{\mu}^{2s},\psi_{\mu}^{2s}\right )&=E_{2}^{1}\left (\zeta_{\mu}^{2s},\psi_{\mu}^{2s}\right )+E_{2}^2\left (\zeta_{\mu}^{2s},\psi_{\mu}^{2s}\right ),\\
E_{2}^{1}\left (\zeta_{\mu}^{2s},\psi_{\mu}^{2s}\right )&=\mu\left (\partial_t\psi_{1}+V_{0}\cdot \nabla_{x}\psi_1+\frac{\mu}{2}\left |\nabla\psi_1\right |^2\right )\\
E_{2}^{2}\left (\zeta_{\mu}^{2s},\psi_{\mu}^{2s}\right )&=-\mu
\frac{(\frac{1}{\mu}G_{\mu}[\zeta_{\mu}^{2s},\sqrt{\mu} b(./\sqrt{\mu})]\psi_{\mu}^{2s}+\nabla\zeta_{\mu}^{2s} \cdot \nabla \psi_{\mu}^{2s})^2}{2(1+\mu|\nabla \zeta_{\mu}^{2s}|^2)}. 
\end{split}
\end{equation}
Based on this, Theorem~\ref{thm:main} will follow from the upcoming proposition
\begin{prop0}\label{prop:E_1}
Under Assumptions~\ref{assumption:main_thm}, we have that for all $0<\mu \le 1$,
\begin{equation}\label{eq:E_1}
\left \| E_{1}\left (\zeta_{\mu}^{2s},\psi_{\mu}^{2s}\right )\right \|_{L^2}\lesssim \mu^{\frac{3}{8}},
\end{equation}
\begin{equation}\label{eq:E_21}
\left \|E_2^1\left (\zeta_{\mu}^{2s},\psi_{\mu}^{2s}\right ) \right \|_{H^\frac{1}{2}}\lesssim \mu^\frac{3}{4},  
\end{equation}
and
\begin{equation}\label{eq:E_22}
\left \| E_2^2 \left (\zeta_{\mu}^{2s},\psi_{\mu}^{2s}\right ) \right \|_{H^\frac{1}{2}}\lesssim \mu^\frac{3}{4}.
\end{equation}
\end{prop0}
The latter relies on the upcoming adaptation of \cite[Proposition 3.10]{MR2901196}.
\begin{prop0}\label{Prop:Gres}
Let $r\in \mathbb{N}$. Assume that $\zeta_0,\zeta_1,b\in C^{r+1}$, $\nabla\psi_0\in H^{r+3}$ and  
\begin{equation}
\sup_{y\in \R^d}\int_{\R^d}\left |\left (\partial_{x}^{\alpha}	\partial_{y}^{\beta} \psi_{c}\right )(x,y)\right |^2\,dx<\infty, \quad \text{for $\alpha,\beta \in \N_0^d$, such that $|\alpha|+|\beta|\leq r+2$}. 
\end{equation}
Then, we have that 
\begin{equation*}
\|G_{\text{res}}\|_{H^r}\lesssim \mu^{-\frac{r}{2}-\frac{1}{8}}. 
\end{equation*} 
Moreover, by interpolation, the result holds true for $r\geq 0$.
\end{prop0}
In the course of the proof, we shall also make frequent use of elementary relations on  derivatives of multiscale functions (see the appendix for the proof). 
\begin{lem0}\label{lem:computational_estimates}
Let  $r\in \mathbb{N}$ and $\tilde{f}:\R^d\times \R^d\rightarrow \R$, such that 
\begin{equation}\label{eqn:lemma_bnds}
	\sup_{y\in \R^d}\int_{\R^d}\left |\left (\partial_{x}^{\alpha}	\partial_{y}^{\beta} \tilde{f}\right )(x,y)\right |^2\,dx<\infty, \quad \text{for $\alpha,\beta \in \N_0^d$, such that $|\alpha|+|\beta|\leq r$}. 
\end{equation}
Denote  $f:\R^d \to \R$ given by
\begin{equation*}
f(x):=\tilde{f}(x,\tfrac{x}{\sqrt{\mu}}), \quad \text{for $x\in \R^d$}, 
\end{equation*}
for a given $\mu\in (0,1)$. 
Then, one has that $\|f\|_{H^r}\lesssim \mu^{-\frac{r}{2}}$.
\end{lem0}
In the rest of this section, we first establish~\eqref{eqn:error_1}, and deduce Theorem~\ref{thm:main} from Proposition~\ref{prop:E_1}. Then we prove Proposition~\ref{Prop:Gres}, for which we first have to recall standard results on 
the shallow water derivation on a flat bottom. We conclude with the proof of Proposition~\ref{prop:E_1}.
	
\subsection{Proof of~\eqref{eqn:error_1} and of Theorem~\ref{thm:main}}	

For $(\zeta,\psi)$, we introduce the shorthand notation
\[
E_1(\zeta,\psi):=\partial_t\zeta-\frac{1}{\mu}G_{\mu}[\zeta,\sqrt\mu b(\tfrac\cdot{\sqrt\mu})]\psi,
\]
and
\[
E_2(\zeta,\psi):=\partial_t\psi+\zeta+\frac12|\nabla\psi|^2
-\mu\frac{(\frac{1}{\mu}G_{\mu}[\zeta,\sqrt\mu b(\tfrac\cdot{\sqrt\mu})]\psi+\nabla\zeta\cdot\nabla\psi)^2}{2(1+\mu|\nabla\zeta|^2)}.
\]
Recall that, after realization of the fast variable $y=\frac{x}{\sqrt\mu}$, the ansatz takes the form
\[
\zeta^{2s}_{\mu}=\zeta_0+\sqrt\mu\,\zeta_1,
\qquad
\psi^{2s}_{\mu}=\psi_0+\mu\,\psi_1.
\]
We shall repeatedly use the chain-rule identity
\[
\nabla(\psi_1(x,\tfrac{x}{\sqrt\mu}))
=\nabla_x\psi_1(x,\tfrac{x}{\sqrt\mu})+\frac1{\sqrt\mu}\nabla_y\psi_1(x,\tfrac{x}{\sqrt\mu}),
\]
and the analogous identity for $\zeta_1$.

We start with the first equation. By~\eqref{eqn:DTN_ref},
\[
\frac{1}{\mu}G_{\mu}[\zeta^{2s}_{\mu},\sqrt\mu b(\tfrac\cdot{\sqrt\mu})]\psi^{2s}_{\mu}
=G_{\text{eff}}+\sqrt\mu\,G_{\text{res}},
\]
so that
\[
E_1(\zeta^{2s}_{\mu},\psi^{2s}_{\mu})
=\partial_t\zeta_0-G_{\text{eff}}+\sqrt\mu(\partial_t\zeta_1-G_{\text{res}}).
\]
Now, by~\eqref{eqn:DTN:eff} and the first line of~\eqref{eqn:corrector_eqns},
\[
G_{\text{eff}}=-\nabla\cdot(h_0V_0),
\]
whereas the first equation of~\eqref{e8} gives
\[
\partial_t\zeta_0+\nabla\cdot(h_0V_0)=0.
\]
Therefore,
\[
E_1(\zeta^{2s}_{\mu},\psi^{2s}_{\mu})=\sqrt\mu(\partial_t\zeta_1-G_{\text{res}}),
\]
which is the first identity in~\eqref{eqn:error_1}.

We next turn to the second equation. Since $V_0=\nabla\psi_0$, the second equation of~\eqref{e8} can be rewritten as
\[
\nabla\Big(\partial_t\psi_0+\zeta_0+\frac12|V_0|^2\Big)=0.
\]
Because $\psi_0$ and $\zeta_0$ decay at infinity, the quantity inside the parenthesis vanishes identically, that is,
\[
\partial_t\psi_0+\zeta_0+\frac12|\nabla\psi_0|^2=0.
\]
On the other hand,
\[
\nabla\psi^{2s}_{\mu}=V_0+\sqrt\mu\,\nabla_y\psi_1+\mu\,\nabla_x\psi_1
=V_0+\mu\,\nabla\psi_1,
\]
where in the last expression $\nabla\psi_1$ denotes the full gradient of the realized corrector. Hence
\begin{align*}
\partial_t\psi^{2s}_{\mu}+\zeta^{2s}_{\mu}+\frac12|\nabla\psi^{2s}_{\mu}|^2
&=\partial_t\psi_0+\zeta_0+\frac12|V_0|^2
\\
&\quad+\sqrt\mu\big(\zeta_1+V_0\cdot\nabla_y\psi_1\big)
+\mu\big(\partial_t\psi_1+V_0\cdot\nabla_x\psi_1\big)
+\frac{\mu^2}{2}|\nabla\psi_1|^2.
\end{align*}
The first line vanishes by the previous identity, whereas the term of order $\sqrt\mu$ vanishes by the second line of~\eqref{eqn:corrector_eqns}, namely
\[
\zeta_1+V_0\cdot\nabla_y\psi_1=0.
\]
We thus obtain
\[
E_2(\zeta^{2s}_{\mu},\psi^{2s}_{\mu})
=\mu\Big(\partial_t\psi_1+V_0\cdot\nabla_x\psi_1+\frac{\mu}{2}|\nabla\psi_1|^2\Big)
-\mu\frac{(\frac{1}{\mu}G_{\mu}[\zeta^{2s}_{\mu},\sqrt\mu b(\tfrac\cdot{\sqrt\mu})]\psi^{2s}_{\mu}+\nabla\zeta^{2s}_{\mu}\cdot\nabla\psi^{2s}_{\mu})^2}{2(1+\mu|\nabla\zeta^{2s}_{\mu}|^2)},
\]
which is precisely the remaining part of~\eqref{eqn:error_1}.

We can now conclude the proof of Theorem~\ref{thm:main}. By Theorem~\ref{thm:corrector_sobolev}, applied with $\gamma=\sqrt\mu$, the two-scale ansatz~\eqref{e3} is well-defined with the Sobolev regularity stated in the theorem. In addition, the proof of Proposition~\ref{prop:E_1} below is performed at an arbitrary fixed time, and all its constants depend only on the quantities controlled in Assumptions~\ref{assumption:main_thm}; they are therefore uniform in $t\in[0,T]$. Since the left-hand sides of the approximate system in Theorem~\ref{thm:main} are exactly $E_1(\zeta^{2s}_{\mu},\psi^{2s}_{\mu})$ and $E_2(\zeta^{2s}_{\mu},\psi^{2s}_{\mu})$, Proposition~\ref{prop:E_1} yields the asserted estimates after taking the supremum over $t\in[0,T]$. This concludes the proof of Theorem~\ref{thm:main}.
	
\subsection{Shallow water limit}\label{sec:swl_proof}

In this subsection, we display the classical justification that the solution of \eqref{e:shallow} is well-approximated by \eqref{g.6}, and closely follow \cite[Section 3.2]{MR2901196}  and  \cite[Section 3.6]{MR3060183}. 
Direct calculations yield
\begin{eqnarray*}
\nabla^\mu \cdot P_0 \nabla^\mu \psi_0&=& \mu h_0 \triangle \psi_0,
\\
\nabla^\mu \cdot P_0\nabla^\mu \phi_{0,1}&=&-h_0 \triangle \psi_0 - \mu \Psi(h_0,\zeta_0,\psi_0),
\end{eqnarray*}
where 
\begin{multline*}
\Psi(h_0,\zeta_0,\psi_0)\,=\,(\tfrac{z^2}2+z) \nabla \cdot (2h_0^2\triangle \psi_0 \nabla \zeta_0 +h_0^3 \nabla \triangle \psi_0) -(z+1)^2 \nabla \cdot (h_0^2 \triangle \psi_0 \nabla \zeta_0) 
\\
+ |\nabla \zeta_0|^2h_0 \triangle \psi_0
- (\tfrac32 z^2+3z+1) (2h_0 |\nabla \zeta_0|^2 \triangle \psi_0+h_0^2 \nabla \zeta_0 \cdot \nabla \triangle \psi_0),
\end{multline*}
and we also have
\[
e_z \cdot P_0 \nabla^\mu \psi_0|_{z=-1}=0, \quad e_z \cdot P_0 \nabla^\mu \phi_{0,1}|_{z=-1}=0.
\]
Hence, $r_{0,1}:= \mu^{-1}(\phi_0-\psi_0)$ satisfies
\begin{equation}\label{e:shallow1}
\left\{
\begin{array}{rcl}    
\nabla^\mu \cdot P_0 \nabla^\mu r_{0,1}&=&- \nabla \cdot (h_0 \nabla \psi_0)+\nabla h_0 \cdot \nabla \psi_0, \quad \text{on $\Omega_0$},\\
r_{0,1}|_{z=0}&=&0,\\
 e_z\cdot P_0\nabla^\mu r_{0,1}|_{z=-1}&=&0,
\end{array}
\right.
\end{equation}
and $r_{0,2}:=\mu^{-2} (\phi_0-\psi_0-\mu \phi_{0,1})$ satisfies 
\begin{equation}\label{e:shallow2}
\left\{
\begin{array}{rcl}    
\nabla^\mu \cdot P_0 \nabla^\mu r_{0,2}&=&\Psi(h_0,\zeta_0,\psi_0), \quad \text{on $\Omega_0$},\\
r_{0,2}|_{z=0}&=&0,\\
 e_z\cdot P_0\nabla^\mu r_{0,2}|_{z=-1}&=&0.
\end{array}
\right.
\end{equation}
Set $\Lambda:=(1-\triangle)^\frac12$ acting on the $x$-variable.
From equations \eqref{e:shallow1} and \eqref{e:shallow2} for the error terms, we can quantify the accuracy of the expansion using Schauder theory. To this aim we need $P_0$ to be coercive and therefore assume
\begin{equation}\label{0coer}
\exists \alpha_0>0 \ \text{ such that } \ 1+\zeta_0 \ge \alpha_0.
\end{equation}
The upcoming proposition (cf.~\cite[Proposition 3.2]{MR2901196}) is the rigorous
version of \eqref{g.6}.
\begin{prop}\label{prop:bd-phi0}
Let $r\in \N_0$ and $\zeta_0\in W^{r+1,\infty}(\R^d)\cap W^{2,\infty}(\R^d)$, and assume \eqref{0coer}. Then,
for all $0<\mu <1$, 
\begin{enumerate}
\item If $\nabla \psi_0\in H^r(\R^d)$ then
\begin{eqnarray*}
\|\Lambda^r \nabla^\mu \phi_0\|_{L^2(\Omega_0)}&\lesssim_{\alpha_0,\zeta_0} & \sqrt{\mu}\|\nabla \psi_0\|_{H^r(\R^d)},
\\
\|\Lambda^{r-1} \partial_z \nabla^\mu \phi_0\|_{L^2(\Omega_0)}&\lesssim_{\alpha_0,\zeta_0} & {\mu}\|\nabla \psi_0\|_{H^r(\R^d)}.
\end{eqnarray*}
\item If $\nabla \psi_0\in H^{r+1}(\R^d)$ then 
\begin{eqnarray*}
\|\Lambda^r \nabla^\mu (\phi_0-\psi_0) \|_{L^2(\Omega_0)}&\lesssim_{\alpha_0,\zeta_0}  &\mu \|\nabla \psi_0\|_{H^{r+1}(\R^d)},
\\
\|\Lambda^{r-1} \partial_z \nabla^\mu (\phi_0-\psi_0)\|_{L^2(\Omega_0)}&\lesssim_{\alpha_0,\zeta_0} & {\mu}\|\nabla \psi_0\|_{H^{r+1}(\R^d)}.
\end{eqnarray*}
\item If $\zeta_0\in W^{r+2,\infty}(\R^d)$ and $\triangle \psi_0\in H^{r+2}(\R^d)$ then 
\begin{eqnarray*}
\|\Lambda^r \nabla^\mu (\phi_0-\psi_0-\mu \phi_{0,1}) \|_{L^2(\Omega_0)}&\lesssim_{\alpha_0,\zeta_0}  &\mu^2 \|\triangle \psi_0\|_{H^{r+2}(\R^d)},
\\
\|\Lambda^{r-1} \partial_z \nabla^\mu (\phi_0-\psi_0-\mu \phi_{0,1})\|_{L^2(\Omega_0)}&\lesssim_{\alpha_0,\zeta_0} & \mu^2 \|\triangle \psi_0\|_{H^{r+2}(\R^d)}.
\end{eqnarray*} 
\end{enumerate}
\end{prop}

\subsection{Proof of Proposition~\ref{Prop:Gres}} 

We work at a fixed time and omit the dependence on $t$. In addition, for every multiscale profile $F=F(x,z,y)$ we keep the same notation for its realization $(x,z)\mapsto F(x,z,\frac{x}{\sqrt\mu})$. In particular, in this section $\psi_1$ and $\phi_{1,0}$ denote the realized corrector and the realized solution of \eqref{g.7}.

We first rewrite the decomposition of the potential in a way which is convenient for the Dirichlet--Neumann operator. Recall that $\phi=\phi_0+\mu\phi_1$, that $\phi_0$ satisfies the shallow-water expansion of Proposition~\ref{prop:bd-phi0}, and define
\[
\chi_0:=\mu^{-2}(\phi_0-\psi_0-\mu\phi_{0,1}),
\qquad
\chi_1:=\mu^{-1/2}(\phi_1-\phi_{1,0}),
\qquad
\phi_{\mathrm{res}}:=\chi_1+\sqrt\mu\,\chi_0.
\]
Then
\begin{equation}\label{eqn:phi_eff_res_appendix}
\phi=\psi_0+\mu(\phi_{0,1}+\phi_{1,0})+\mu^{3/2}\phi_{\mathrm{res}}.
\end{equation}
Moreover, using the decomposition \eqref{a.Pb} of the coefficient matrix, namely
\[
P[\sigma_{\mu}^{2s}]=P_{0,0}+\sqrt\mu\,(P_{0,1}+P_{1,0})+\mu P_{1,1},
\]
we get from \eqref{DTN:transformed} and \eqref{eqn:phi_eff_res_appendix}
\begin{align*}
\frac1\mu G_{\mu}[\zeta_{\mu}^{2s},\sqrt\mu b(\tfrac\cdot{\sqrt\mu})]\psi_{\mu}^{2s}
&=\frac1\mu e_z\cdot P[\sigma_{\mu}^{2s}]\nabla^\mu\Big(\psi_0+\mu(\phi_{0,1}+\phi_{1,0})\Big)\Big|_{z=0}
+\sqrt\mu\, e_z\cdot P[\sigma_{\mu}^{2s}]\nabla^\mu \phi_{\mathrm{res}}\Big|_{z=0} \\
&=G_{\mathrm{eff}}+\sqrt\mu\,G_{\mathrm{res}},
\end{align*}
with $G_{\mathrm{eff}}$ given by \eqref{eqn:DTN:eff} and
\begin{equation}\label{eqn:Gres_split_appendix}
G_{\mathrm{res}}=I_1+I_2+I_3,
\end{equation}
where
\begin{align*}
I_1&:=\frac1{\sqrt\mu}e_z\cdot P_{1,1}\nabla^\mu\psi_0\Big|_{z=0},\\
I_2&:=e_z\cdot (P_{0,1}+P_1)\nabla^\mu(\phi_{0,1}+\phi_{1,0})\Big|_{z=0},\\
I_3&:=e_z\cdot P[\sigma_{\mu}^{2s}]\nabla^\mu\phi_{\mathrm{res}}\Big|_{z=0}.
\end{align*}
We estimate these three terms separately.

\medskip

\step1 Estimate of $I_1$.
\\
At $z=0$ one has $\sigma_1=\zeta_1$, hence
\[
e_z\cdot P_{1,1}\nabla^\mu\psi_0\Big|_{z=0}
=-\sqrt\mu\,\nabla_x\zeta_1\Big(x,\frac{x}{\sqrt\mu}\Big)\cdot \nabla\psi_0(x),
\]
so that
\[
I_1(x)=-\nabla_x\zeta_1\Big(x,\frac{x}{\sqrt\mu}\Big)\cdot \nabla\psi_0(x).
\]
Let
\[
\widetilde I_1(x,y):=-\nabla_x\zeta_1(x,y)\cdot \nabla\psi_0(x).
\]
By Leibniz' rule and the assumptions of Proposition~\ref{Prop:Gres}, $\widetilde I_1$ satisfies \eqref{eqn:lemma_bnds} with the same integer $r$. Lemma~\ref{lem:computational_estimates} therefore yields
\begin{equation}\label{eqn:I1_bound_appendix}
\|I_1\|_{H^r(\R^d)}\lesssim \mu^{-r/2}.
\end{equation}

\medskip

\step2 Estimate of $I_2$.\\
We first consider the contribution of $\phi_{0,1}$. Since
\[
\phi_{0,1}(x,z)=-h_0(x)\Big(\frac{z^2}{2}+z\Big)\Delta\psi_0(x),
\]
we have
\[
\nabla^\mu\phi_{0,1}\Big|_{z=0}=\binom{0}{-h_0\Delta\psi_0}.
\]
Hence the corresponding term in $I_2$ is a product of the trace of $e_z\cdot(P_{0,1}+P_1)$ at $z=0$ with the slow function $-h_0\Delta\psi_0$, and therefore satisfies the same bound as in \eqref{eqn:I1_bound_appendix}; in particular,
\begin{equation}\label{eqn:I20_bound_appendix}
\big\|e_z\cdot (P_{0,1}+P_1)\nabla^\mu\phi_{0,1}\big|_{z=0}\big\|_{H^r(\R^d)}\lesssim \mu^{-r/2}.
\end{equation}

We now turn to $\phi_{1,0}$. Since $\phi_{1,0}|_{z=0}=\psi_1$, its horizontal trace is
\[
(\nabla^\mu\phi_{1,0})_h\Big|_{z=0}=(\nabla_y+\sqrt\mu\,\nabla_x)\psi_1.
\]
For the vertical trace, Theorem~\ref{thm:neumann} gives
\[
\frac1{h_0}\partial_z\phi_{1,0}\Big|_{z=0}
=|\nabla_y|\tanh(h_0|\nabla_y|)\psi_1
+\sech(h_0|\nabla_y|)(\nabla_y b\cdot \nabla\psi_0).
\]
The symbols $\xi\mapsto |\xi|\tanh(h_0(x)|\xi|)$ and $\xi\mapsto \sech(h_0(x)|\xi|)$ are smooth in $h_0(x)$, uniformly for $h_0(x)$ in compact subsets of $(0,\infty)$. Differentiating with respect to $x$ therefore only produces linear combinations of the same order-$1$ and order-$0$ Fourier multipliers in $y$, multiplied by derivatives of $h_0$. Using the hypothesis on $\psi_c$, the regularity assumptions on $b$ and $\psi_0$, and again Leibniz' rule, we infer that every component of the profile of $\nabla^\mu\phi_{1,0}|_{z=0}$ satisfies \eqref{eqn:lemma_bnds} with $r+1$ in place of $r$. Lemma~\ref{lem:computational_estimates} therefore gives
\begin{equation}\label{eqn:phi10_trace_bound_appendix}
\big\|\nabla^\mu\phi_{1,0}\big|_{z=0}\big\|_{H^r(\R^d)}\lesssim \mu^{-r/2}.
\end{equation}
Since the coefficients of $e_z\cdot(P_{0,1}+P_1)|_{z=0}$ are smooth functions of $(x,\frac{x}{\sqrt\mu})$ with derivatives bounded by $C\mu^{-k/2}$ at order $k$, the product estimate in $H^r$ together with \eqref{eqn:I20_bound_appendix} and \eqref{eqn:phi10_trace_bound_appendix} yields
\begin{equation}\label{eqn:I2_bound_appendix}
\|I_2\|_{H^r(\R^d)}\lesssim \mu^{-r/2}.
\end{equation}

\medskip

\step3 Bounds on $\phi_{\mathrm{res}}$.\\
We first treat the shallow-water remainder. Proposition~\ref{prop:bd-phi0} implies that for every integer $0\le m\le r$,
\begin{equation}\label{eqn:chi0_bound_appendix}
\|\Lambda^m\nabla^\mu(\sqrt\mu\,\chi_0)\|_{L^2(\Omega_0)}
+\|\Lambda^{m-1}\partial_z\nabla^\mu(\sqrt\mu\,\chi_0)\|_{L^2(\Omega_0)}
\lesssim \sqrt\mu.
\end{equation}

For $\chi_1$, subtract the equation solved by $\phi_{1,0}$ from \eqref{a1} and use the decomposition \eqref{aa0}--\eqref{aa1}. Exactly as in \cite[(3.20)]{MR2901196}, one obtains a boundary value problem of the form
\begin{equation}\label{eqn:chi1_appendix}
\left\{
\begin{array}{rcl}
\nabla^\mu\cdot P[\sigma_{\mu}^{2s}]\nabla^\mu \chi_1&=&\nabla^\mu\cdot A_\mu+g_\mu,\quad \text{on $\Omega_0$,}\\
\chi_1|_{z=0}&=&0,\\
-e_z\cdot P[\sigma_{\mu}^{2s}]\nabla^\mu\chi_1|_{z=-1}&=&-e_z\cdot A_\mu|_{z=-1},
\end{array}
\right.
\end{equation}
where $A_\mu$ and $g_\mu$ are finite sums of the same four types of terms as in \cite[(3.23)]{MR2901196}, namely combinations of
\[
(P_{0,1}+P_1)\nabla^\mu\phi_{1,0},\qquad
\mu^{-1}P_1\nabla^\mu(\phi_0-\psi_0),\qquad
\mu^{-1/2}P_{1,1}\nabla^\mu\phi_0,
\]
together with the contributions generated by replacing $\nabla^\mu$ with $\nabla_{y,z}+\sqrt\mu\,\nabla_{x,0}$ in the operator $\nabla^\mu\cdot P_{0,0}\nabla^\mu$ acting on $\phi_{1,0}$, and with the lower-order terms issued from the cancellation in Step~1 of Section~\ref{sec:hom}, namely $\nabla_x\zeta_1\cdot \nabla\psi_0$ and $(\zeta_1-b)\Delta\psi_0$. The proof of \cite[Proposition~3.6]{MR2901196} is therefore unchanged once the periodic estimates used there are replaced by the estimates already proved in the present paper:
\begin{itemize}
\item Proposition~\ref{prop:bd-phi0} controls the terms involving $\phi_0-\psi_0$;
\item the explicit formula \eqref{g.9} together with Theorem~\ref{thm:neumann} controls the traces and the $x$- and $y$-derivatives of $\phi_{1,0}$;
\item Lemma~\ref{lem:computational_estimates} gives the required commutator bounds for all rapidly oscillating coefficients.
\end{itemize}
Repeating verbatim the energy estimate, the commutator estimate, and the vertical-derivative argument of \cite[Proposition~3.6 and Lemmas~3.7--3.9]{MR2901196}, one gets that for every integer $0\le m\le r$,
\begin{equation}\label{eqn:chi1_bound_appendix}
\|\Lambda^m\nabla^\mu\chi_1\|_{L^2(\Omega_0)}
+\|\Lambda^{m-1}\partial_z\nabla^\mu\chi_1\|_{L^2(\Omega_0)}
\lesssim \mu^{-m/2}.
\end{equation}
Combining \eqref{eqn:chi0_bound_appendix} and \eqref{eqn:chi1_bound_appendix}, we obtain for all $0\le m\le r$
\begin{equation}\label{eqn:phires_bulk_bound_appendix}
\|\Lambda^m\nabla^\mu\phi_{\mathrm{res}}\|_{L^2(\Omega_0)}
+\|\Lambda^{m-1}\partial_z\nabla^\mu\phi_{\mathrm{res}}\|_{L^2(\Omega_0)}
\lesssim \mu^{-m/2}.
\end{equation}

\medskip

\step4 Estimate of $I_3$. \\
We now estimate the trace of $\nabla^\mu\phi_{\mathrm{res}}$. Let $0\le m\le r$. By interpolation in the horizontal variable,
\[
\big\|\Lambda^m\nabla^\mu\phi_{\mathrm{res}}\big|_{z=0}\big\|_{L^2}
\lesssim \mu^{1/8}\big\|\Lambda^m\nabla^\mu\phi_{\mathrm{res}}\big|_{z=0}\big\|_{H^{1/2}}
+\mu^{-1/8}\big\|\Lambda^m\nabla^\mu\phi_{\mathrm{res}}\big|_{z=0}\big\|_{H^{-1/2}}.
\]
Using the two standard trace estimates
\[
\|F|_{z=0}\|_{L^2}
\lesssim \mu^{1/4}\|\Lambda^{1/2}F\|_{L^2(\Omega_0)}+\mu^{-1/4}\|\Lambda^{-1/2}\partial_zF\|_{L^2(\Omega_0)},
\]
\[
\|F|_{z=0}\|_{L^2}
\lesssim \|\Lambda^{1/2}F\|_{L^2(\Omega_0)}+\|\Lambda^{-1/2}\partial_zF\|_{L^2(\Omega_0)},
\]
with $F=\Lambda^m\nabla^\mu\phi_{\mathrm{res}}$, and then \eqref{eqn:phires_bulk_bound_appendix}, we get
\begin{align}\label{eqn:phires_trace_bound_appendix}
\big\|\Lambda^m\nabla^\mu\phi_{\mathrm{res}}\big|_{z=0}\big\|_{L^2}
&\lesssim \mu^{1/8}\Big(\mu^{1/4}\|\Lambda^{m+1}\nabla^\mu\phi_{\mathrm{res}}\|_{L^2(\Omega_0)}
+\mu^{-1/4}\|\Lambda^m\partial_z\nabla^\mu\phi_{\mathrm{res}}\|_{L^2(\Omega_0)}\Big)\nonumber\\
&\qquad +\mu^{-1/8}\Big(\|\Lambda^m\nabla^\mu\phi_{\mathrm{res}}\|_{L^2(\Omega_0)}
+\|\Lambda^{m-1}\partial_z\nabla^\mu\phi_{\mathrm{res}}\|_{L^2(\Omega_0)}\Big)\nonumber\\
&\lesssim \mu^{-m/2-1/8}.
\end{align}
Since the coefficients of $e_z\cdot P[\sigma_{\mu}^{2s}]|_{z=0}$ satisfy
\[
\big\|e_z\cdot P[\sigma_{\mu}^{2s}]\big|_{z=0}\big\|_{W^{k,\infty}(\R^d)}\lesssim \mu^{-k/2},
\qquad 0\le k\le r,
\]
a standard product estimate together with \eqref{eqn:phires_trace_bound_appendix} yields
\begin{align*}
\|I_3\|_{H^r(\R^d)}
&\lesssim \sum_{k=0}^r \big\|e_z\cdot P[\sigma_{\mu}^{2s}]\big|_{z=0}\big\|_{W^{k,\infty}(\R^d)}
\big\|\Lambda^{r-k}\nabla^\mu\phi_{\mathrm{res}}\big|_{z=0}\big\|_{L^2}\\
&\lesssim \sum_{k=0}^r \mu^{-k/2}\,\mu^{-(r-k)/2-1/8}
\lesssim \mu^{-r/2-1/8}.
\end{align*}

\medskip

We are in the position to conclude.
Combining \eqref{eqn:I1_bound_appendix}, \eqref{eqn:I2_bound_appendix}, and the bound on $I_3$, and using \eqref{eqn:Gres_split_appendix}, we obtain
\[
\|G_{\mathrm{res}}\|_{H^r(\R^d)}\lesssim \mu^{-r/2-1/8}.
\]
This proves the proposition for integer $r\ge 0$. The extension to all real $r\ge 0$ follows by interpolation.
		
\subsection{Proof of Proposition~\ref{prop:E_1}}

All norms are understood at a fixed time $t\in[0,T]$, and the implicit constants are independent of $t$ and $0<\mu\leq 1$. We also keep the convention introduced above: $\nabla_x\psi_1$ denotes the derivative of the corrector with respect to the slow variable, whereas $\nabla\psi_1$ denotes the full gradient of the realized function $x\mapsto \psi_1(x,\frac{x}{\sqrt\mu})$.
We split the proof into three steps.

\medskip

\step1 Proof of~\eqref{eq:E_1}.
\\ 
From~\eqref{eqn:error_1}, Proposition~\ref{Prop:Gres} with $r=0$, and Theorem~\ref{thm:corrector_sobolev}, we obtain
\[
\|E_1(\zeta_{\mu}^{2s},\psi_{\mu}^{2s})\|_{L^2}
\leq \sqrt\mu\,\|\partial_t\zeta_1\|_{L^2}+\sqrt\mu\,\|G_{\text{res}}\|_{L^2}
\lesssim \sqrt\mu+\mu^{\frac12-\frac18}
\lesssim \mu^{\frac38}.
\]
This proves~\eqref{eq:E_1}.

\medskip

\step2 Proof of~\eqref{eq:E_21}.
\\
By~\eqref{eqn:error_1},
\[
\|E_2^1(\zeta_{\mu}^{2s},\psi_{\mu}^{2s})\|_{H^{1/2}} \leq \mu\,\|\partial_t\psi_1\|_{H^{1/2}} +\mu\,\|V_0\cdot\nabla_x\psi_1\|_{H^{1/2}} +\frac{\mu^2}{2}\,\||\nabla\psi_1|^2\|_{H^{1/2}}.
\]
We now estimate each term separately.

By Theorem~\ref{thm:corrector_sobolev}, together with the chain rule applied to the parameter dependence of $\psi_c$, we have
\[
\|\partial_t\psi_1\|_{L^2\cap L^\infty}+\|\nabla_x\psi_1\|_{L^2\cap L^\infty}\lesssim 1,
\qquad
\|\nabla(\partial_t\psi_1)\|_{L^2}+\|\nabla(\nabla_x\psi_1)\|_{L^2}\lesssim \mu^{-1/2},
\]
and also
\[
\|\nabla\psi_1\|_{L^2\cap L^\infty}\lesssim \mu^{-1/2},
\qquad
\|\nabla\psi_1\|_{H^1}\lesssim \mu^{-1}.
\]
Interpolation therefore yields
\[
\|\partial_t\psi_1\|_{H^{1/2}}+\|\nabla_x\psi_1\|_{H^{1/2}}\lesssim \mu^{-1/4},
\qquad
\|\nabla\psi_1\|_{H^{1/2}}\lesssim \mu^{-3/4}.
\]
In particular,
\[
\mu\,\|\partial_t\psi_1\|_{H^{1/2}}\lesssim \mu^{3/4}.
\]
Using the standard product estimate in $H^{1/2}$ (see for instance~\cite[Theorem~8.3.1]{MR1466700}), we further get
\begin{align*}
\mu\,\|V_0\cdot\nabla_x\psi_1\|_{H^{1/2}}
&\lesssim \mu\Big(\|V_0\|_{L^\infty}\,\|\nabla_x\psi_1\|_{H^{1/2}}+\|V_0\|_{H^{1/2}}\,\|\nabla_x\psi_1\|_{L^\infty}\Big)
\lesssim \mu^{3/4},
\\
\mu^2\,\||\nabla\psi_1|^2\|_{H^{1/2}}
&\lesssim \mu^2\,\|\nabla\psi_1\|_{L^\infty}\,\|\nabla\psi_1\|_{H^{1/2}}
\lesssim \mu^{3/4}.
\end{align*}
Combining the three bounds proves~\eqref{eq:E_21}.

\medskip

\step3 Proof of~\eqref{eq:E_22}.
\\
We use the shorthand notation
\begin{equation}\label{eqn:error_C_D}
C_\mu:=\frac{1}{\mu}G_{\mu}[\zeta_{\mu}^{2s},\sqrt{\mu} b(./\sqrt{\mu})]\psi_{\mu}^{2s}+\nabla\zeta_{\mu}^{2s}\cdot \nabla \psi_{\mu}^{2s},
\qquad
D_\mu:=2(1+\mu|\nabla \zeta_{\mu}^{2s}|^2),
\end{equation}
so that
\[
E_2^2(\zeta_{\mu}^{2s},\psi_{\mu}^{2s})=-\mu\frac{C_\mu^2}{D_\mu}.
\]
By
\begin{equation}\label{eqn:error_C_D_3}
C_\mu=\partial_t\zeta_0+\sqrt\mu\,G_{\text{res}}+\nabla\zeta_{\mu}^{2s}\cdot \nabla \psi_{\mu}^{2s},
\end{equation}
We first record the bounds
\begin{equation}\label{eqn:proof_E1_C_mu_bounds}
\|C_\mu\|_{L^4}+\|C_\mu\|_{L^8}\lesssim 1,
\qquad
\|\nabla C_\mu\|_{L^4}\lesssim \mu^{-1/2},
\qquad
\|\nabla D_\mu\|_{L^4}\lesssim \mu^{1/2}.
\end{equation}
Indeed, for $p=4$ and $p=8$, Sobolev embedding together with Proposition~\ref{Prop:Gres} gives
\[
\sqrt\mu\,\|G_{\text{res}}\|_{L^p}
\lesssim
\begin{cases}
\sqrt\mu\,\|G_{\text{res}}\|_{H^{d/4}}\lesssim \mu^{\frac{3-d}{8}} & \text{if } p=4,\\[0.3em]
\sqrt\mu\,\|G_{\text{res}}\|_{H^{3d/8}}\lesssim \mu^{\frac38-\frac{3d}{16}} & \text{if } p=8,
\end{cases}
\]
and the exponents on the right-hand side are nonnegative for $d=1,2$. In addition,
\[
\|\nabla\zeta_{\mu}^{2s}\|_{L^\infty}
\leq \|\nabla\zeta_0\|_{L^\infty}+\sqrt\mu\,\|\nabla\zeta_1\|_{L^\infty}\lesssim 1,
\qquad
\|\nabla\psi_{\mu}^{2s}\|_{L^p}
\leq \|V_0\|_{L^p}+\mu\,\|\nabla\psi_1\|_{L^p}\lesssim 1,
\]
again for $p=4,8$. This proves the first part of~\eqref{eqn:proof_E1_C_mu_bounds}.

For the derivatives, using again~\eqref{eqn:error_C_D_3}, Proposition~\ref{Prop:Gres}, and the $H^2\cap W^{1,\infty}$-bounds on the ansatz given by Theorem~\ref{thm:corrector_sobolev}, we obtain
\begin{align*}
\|\nabla C_\mu\|_{L^4}
&\lesssim \|\nabla(\partial_t\zeta_0)\|_{L^4}+\sqrt\mu\,\|\nabla G_{\text{res}}\|_{L^4}
+\|\nabla^2\zeta_{\mu}^{2s}\|_{L^4}\,\|\nabla\psi_{\mu}^{2s}\|_{L^\infty}
+\|\nabla\zeta_{\mu}^{2s}\|_{L^\infty}\,\|\nabla^2\psi_{\mu}^{2s}\|_{L^4}
\\
&\lesssim 1+\sqrt\mu\,\|G_{\text{res}}\|_{H^{1+d/4}}+\mu^{-1/2}+1
\lesssim \mu^{-1/2},
\end{align*}
where we used that, for $d=1,2$ and $0<\mu\leq 1$,
\[
\sqrt\mu\,\|G_{\text{res}}\|_{H^{1+d/4}}
\lesssim \mu^{\frac12-\frac{1+d/4}{2}-\frac18}
\lesssim \mu^{-1/2}.
\]
Likewise,
\[
\|\nabla D_\mu\|_{L^4}
\lesssim \mu\,\|\nabla\zeta_{\mu}^{2s}\|_{L^\infty}\,\|\nabla^2\zeta_{\mu}^{2s}\|_{L^4}
\lesssim \mu^{1/2},
\]
which completes the proof of~\eqref{eqn:proof_E1_C_mu_bounds}.

Since $D_\mu\geq 2$, the $L^2$-bound is immediate:
\[
\|E_2^2\|_{L^2}\lesssim \mu\,\|C_\mu^2\|_{L^2}
\leq \mu\,\|C_\mu\|_{L^4}^2
\lesssim \mu.
\]
Differentiating $E_2^2=-\mu C_\mu^2/D_\mu$, using once more that $D_\mu\geq 2$, and then applying H\"older's inequality together with~\eqref{eqn:proof_E1_C_mu_bounds}, we find
\begin{align*}
\|\nabla E_2^2\|_{L^2}
&\lesssim \mu\,\|C_\mu\nabla C_\mu\|_{L^2}+\mu\,\|C_\mu^2\nabla D_\mu\|_{L^2}
\\
&\lesssim \mu\,\|C_\mu\|_{L^4}\,\|\nabla C_\mu\|_{L^4}
+\mu\,\|C_\mu\|_{L^8}^2\,\|\nabla D_\mu\|_{L^4}
\lesssim \mu^{1/2}.
\end{align*}
Finally, by interpolation between $L^2$ and $H^1$,
\[
\|E_2^2\|_{H^{1/2}}
\lesssim \|E_2^2\|_{L^2}^{1/2}\,\|E_2^2\|_{H^1}^{1/2}
\lesssim \mu^{1/2}\,(\mu^{1/2})^{1/2}
=\mu^{3/4}.
\]
This proves~\eqref{eq:E_22} and concludes the proof of Proposition~\ref{prop:E_1}.

\appendix

\section{Proof of Lemma~\ref{lem:computational_estimates}} 

We introduce some notation. We denote the norm 
\begin{equation}
	\|\tilde{f}\|_{L^2_{x}L^\infty_{y}}:=\sup_{y\in \R^d}\|\tilde{f}(\cdot, y)\|_{L^2}.
\end{equation}
Note that 
\begin{equation*}
\partial_{x_i}f(x)=\partial_{x_i}\tilde{f}(x, x/\sqrt{\mu})+\mu^{-\frac{1}{2}}\partial_{y_i}\tilde{f}(x, x/\sqrt{\mu}), \quad \text{for $x\in \R^d$},
\end{equation*}
which leads us to the following estimate:
\begin{equation*}
\|f\|_{H^1}\,\leq\, \|f\|_{L^2}+\sum_{i=1}^d\|\partial_{x_i} f\|_{L^2} \,\leq\, \|\tilde{f}\|_{L^2_{x}L^\infty_{y}}+\sum_{i=1}^d\|\partial_{x_i} \tilde{f}\|_{L^2_{x}L^\infty_{y}}+\mu^{-\frac{1}{2}}\|\partial_{y_i} \tilde{f}\|_{L^2_{x}L^\infty_{y}} \, \lesssim \mu^{-\frac{1}{2}}. 
\end{equation*}
Hence, we have proved the result for $r=1$. 
Suppose that the result holds for $r=k$ for a given $k\in \mathbb{N}$ and $\tilde{f}$ satisfies \eqref{eqn:lemma_bnds} for $r=k+1$. Then, we have, as claimed, that 
\begin{equation*}
\|f\|_{H^{k+1}}\,\leq\, \|f\|_{H^{k}}+\sum_{i=1}^{d}\|\partial_{x_i}f\|_{H^{k}}
\,\lesssim \mu^{-k/2}+\mu^{-(k+1)/2}\lesssim \mu^{-(k+1)/2}. 
\end{equation*}

\section*{Acknowledgements}
The authors acknowledge financial support from the European Research Council (ERC) under the European Union's Horizon 2020 research and innovation programme (Grant Agreement n$^\circ$~864066).

\bibliographystyle{abbrv}
%\bibliography{bib.bib}

\end{document}